\documentclass[11pt]{amsart}
\usepackage{amssymb}
\usepackage{graphicx}                          

\newtheorem{lemma}{Lemma}[section]
\newtheorem{theorem}{Theorem}[section]

\newenvironment{adfenumerate}{
\begin{enumerate}
\setlength{\itemsep}{0.5mm}
\setlength{\parskip}{0mm}
\setlength{\parsep}{0mm}
}{
\end{enumerate}
}

\newcommand{\adfmod}[1]{~(\mathrm{mod}~#1)}    
\newcommand{\adfPENT}{\mathop{\mathrm{PENT}} } 
\begin{document}
\title{Generalized pentagonal geometries}
\author{A. D. Forbes}
\address{ORCID: 0000-0003-3805-7056}
\email{anthony.d.forbes@gmail.com}
\author{C. G. Rutherford}
\address{LSBU Business School,
London South Bank University,
103 Borough Road,
London SE1 0AA, UK}
\email{c.g.rutherford@lsbu.ac.uk}
\date{\today}

\maketitle
\begin{abstract}
A pentagonal geometry PENT($k$, $r$) is a partial linear space, where
every line is incident with $k$ points,
every point is incident with $r$ lines, and
for each point $x$, there is a line incident with precisely those points that are not collinear with $x$.

Here we generalize the concept by allowing the points not collinear with $x$ to form the point set of
a Steiner system $S(2,k,w)$ whose blocks are lines of the geometry.
\end{abstract}


\section{Introduction}\label{sec:Introduction}
A {\em Steiner system} $S(2,k,w)$ is an ordered pair $(W, \mathcal{B})$ such that
(i) $W$ is a set of {points} with $|W| = w$,
(ii) $\mathcal{B}$ is a set of $k$-subsets of $B$, called {blocks}, and
(iii) each pair $\{x,y\} \subseteq W$, $x \neq y$ is a subset of precisely one block of $\mathcal{B}$.

For the purpose of this paper, a {\em generalized pentagonal geometry} is a uniform, regular partial linear space such that
for each point $x$, there is a Steiner system $(W_x, \mathcal{B}_x)$ where
$W_x$ is the set of points which are not collinear with $x$ and
the blocks of $\mathcal{B}_x$ are lines of the geometry.
By uniformity and regularity the number of points that are collinear with $x$ does not depend on $x$.
Therefore $|W_x|$ must also be constant.

If $k$ is the number of points incident with a given line,
$r$ is the number of lines incident with a given point
and $w$ is the point set size of the Steiner systems, then $k$, $r$ and $w$ are constants,
and we denote a generalized pentagonal geometry with these parameters by $\adfPENT(k,r,w)$.
The Steiner system $S(2,k,w)$ consisting of the points that are not collinear with point $x$ is
called the {\em opposite design} to $x$ and is denoted by $x^\mathrm{opp}$.
The notation $x^\mathrm{opp}$ may refer to the $S(2,k,w)$ itself, or its point set, or its block set,
whichever is appropriate from the context.

To avoid excessive trivialities we always assume that $k \ge 2$.
Also we recognise three simplistic cases: for all $k \ge 2$, there exist Steiner systems $S(2, k, 0)$, $S(2, k, 1)$ and $S(2, k, k)$,
with block set sizes 0, 0 and 1 respectively.
If a Steiner system $S(2,k,w)$ exists, then
the number of blocks is $w(w - 1)/(k(k - 1))$ and the number of blocks that contain a given point is a constant, $(w - 1)/(k - 1)$.
An $S(2,2,w)$ is essentially a complete graph $K_w$, and
an $S(2,3,w)$ is also called a {\em Steiner triple system} of order $w$, STS$(w)$.
An {\em affine plane} of order $n$ is a Steiner system $S(2,n,n^2)$ and
a {\em projective plane} of order $n$ is a Steiner system $S(2, n+1, n^2 + n + 1)$.
We also refer to Steiner systems generally as designs.
A necessary and sufficient condition for the existence of a Steiner system $S(2,k,w)$ with $w \ge 1$ is
$w \equiv 1 \textrm{~or~} k \adfmod{k(k - 1)}$
when $k = 2$ (trivially), $k = 3$, \cite{Kirkman1847}, and $k \in \{4,5\}$, \cite{Hanani1961}.
An affine plane and a projective plane of order $n$ exist whenever $n$ is a prime power.

A $\adfPENT(k,r,0)$ is essentially a Steiner system $S(2, k, (k - 1)r + 1)$, and
a $\adfPENT(k,r,1)$ consists of the points and blocks of a $k$-GDD of type $2^{r(k - 1)/2 + 1}$
as defined in Section~\ref{sec:Generalized pentagonal geometries}.
If $w = k$, the structure is known as a {\em pentagonal geometry} and may be denoted by $\adfPENT(k,r)$.
We must have $r \ge (w-1)/(k-1)$ and, as previously stated, we always assume $k \ge 2$.
For the remainder of the paper, we assume $w \ge k$ in order to avoid $S(2,k,w)$ designs with empty block sets, and
wherever necessary we tacitly assume that the relevant Steiner systems exist.

For convenience, we shall regard lines as sets of points---thus two points are collinear if they are on (i.e.\ elements of) the same line.
Lines are also called blocks, and we refer to the parameter $k$ of a $\adfPENT(k,r,w)$ as the block size.

The {\em deficiency graph} of a generalized pentagonal geometry $\adfPENT(k,r,w)$
has as its vertices the points of the geometry, and
there is an edge $x \sim y$ precisely when $x$ and $y$ are not collinear.
It is clear that the graph must be $w$-regular and triangle-free.

For brevity, we adopt the notation {\em $(w,g)$-graph} for a graph that is $w$-regular and has girth $g$,
and {\em $(w,g^+)$-graph} for a graph that is $w$-regular and has girth at least $g$.
Thus we can say that a $\adfPENT(k,r,w)$ has a deficiency $(w,4^+)$-graph.

The concept of a {pentagonal geometry} was introduced in \cite{BallBambergDevillersStokes2013} to provide a generalization
of the pentagon analogous to the generalization of the polygon as described in \cite{Tits1959} and \cite{FeitHigman1964}.
The geometry described by Ball, Bamberg, Devillers \& Stokes in \cite{BallBambergDevillersStokes2013}
is based on the observation that for each vertex $x$ of a pentagon, the two vertices that are not collinear with $x$ form a line.
According to their definition the pentagon is a $\adfPENT(2,2)$.

Our generalization is based on the observation that a line containing $k$ points
is the single block of the Steiner system $S(2,k,k)$.
It seems natural, therefore, to extend the definition of $x^\mathrm{opp}$ to a general Steiner system $S(2,k,w)$
since these designs have the same relevant property as a single line: any two points are collinear.
According to our definition the pentagon is a $\adfPENT(2,2,2)$.
For a more substantial example, let $H$ be the Hoffman--Singleton graph,
the unique 50-vertex $(7,5)$-graph, \cite{HoffmanSingleton1960}.
Then we know that the neighbours of the vertices of $H$ form the lines of a $\adfPENT(7,7)$, \cite[Section 3]{BallBambergDevillersStokes2013}.
For each vertex $x$ of $H$, construct a Steiner triple system of order 7 on the set of neighbours of $x$.
Since the neighbourhoods of two distinct vertices of $H$ cannot have more than one
common vertex---otherwise the graph would not have girth 5---the block sets of the fifty STS(7)s are mutually disjoint.
The 350 blocks so obtained are the lines of a $\adfPENT(3,21,7)$.
We are also motivated by the limitations of our knowledge concerning pentagonal geometries $\adfPENT(k,r)$ with $k > 7$.
In Section~\ref{sec:Block size 3} we see that it is not especially difficult to find examples of
$\adfPENT(3,r,w)$ with connected deficiency $(w,5^+)$-graphs for $w \in \{9, 13, 15, 19, 21\}$
but we are not aware of any examples of the corresponding $\adfPENT(w,r)$.

If $w \ge k \ge 2$ and there exists an $S(2,k,w)$, then there exists a $\adfPENT(k, (w - 1)/(k - 1), w)$.
In this generalized pentagonal geometry, which by analogy with \cite{BallBambergDevillersStokes2013} we regard as degenerate,
the lines are the blocks of two Steiner systems $S(2,k,w)$ with disjoint point sets and
the deficiency graph is a complete bipartite graph $K_{w,w}$.
Thus for any point $x$ in either one of the Steiner systems, $x^\mathrm{opp}$ is the other system.
The lines of a degenerate geometry can also occur in the line set of a larger generalized pentagonal geometry, and,
again by analogy with \cite{BallBambergDevillersStokes2013}, we refer to the substructure as an {\em opposite design pair},
or when $w = k$, an {\em opposite line pair}.

For the general theory of pentagonal geometries as well as further background material, we refer the reader to \cite{BallBambergDevillersStokes2013}.
The subject was further developed in \cite{GriggsStokes2016}, \cite{ForbesGriggsStokes2020} and \cite{Forbes2020P}.
The existence spectrum was settled for $\adfPENT(2,r)$ in \cite[Corollary 2.6]{BallBambergDevillersStokes2013},
and for $\adfPENT(3,r)$ as follows:\
unconditionally, i.e.\ with no restriction on the number of opposite line pairs, \cite[Theorem 9]{GriggsStokes2016};
for geometries without opposite line pairs, \cite{ForbesGriggsStokes2020};
and for geometries with a given number of opposite line pairs, \cite{ForbesGriggsStokes2020}.
Similar results have been established for $\adfPENT(4,r)$:\
unconditionally, with a few possible exceptions, \cite{ForbesGriggsStokes2020};
for geometries without opposite line pairs, with a few possible exceptions, \cite{Forbes2020P};
and for a given number of opposite line pairs and sufficiently large $r$, \cite{Forbes2020P}.
Reference \cite{Forbes2020P} also deals with block size 5,
where some progress is made towards the solution of the existence spectrum problem for
$\adfPENT(5,r)$ without opposite line pairs.
If opposite line pairs are allowed, the spectrum for $\adfPENT(5,5t + 1)$ was determined up to a few possible exceptions in \cite[Theorem 12]{GriggsStokes2016}.

The situation regarding pentagonal geometries $\adfPENT(k,r)$ with {\em connected} deficiency graphs is far from resolved for $k \ge 3$.
All of the construction methods employed in \cite{GriggsStokes2016}, \cite{ForbesGriggsStokes2020} and \cite{Forbes2020P}
use group divisible designs, as in Theorem~\ref{thm:GDD-basic}, below, and inevitably produce
geometries with deficiency graphs consisting of at least $k$ connected components.
As far as we are aware, the known $\adfPENT(k,r)$ with $k \ge 3$ and connected deficiency graph of girth at least 5
amount to the following:
\begin{enumerate}
\item[(i)] a small number of $\adfPENT(3,r)$ constructed by hand, including
  the Desargues configuration, $\adfPENT(3,3)$, \cite{BallBambergDevillersStokes2013};
\item[(ii)] $\adfPENT(3,r)$ that can be created by hill climbing for $r$ not too large, \cite{ForbesGriggsStokes2020};
\item[(iii)] $\adfPENT(3,r)$ for all $r \equiv 3 \adfmod{6}$, $r \ge 33$, \cite{Forbes2020P};
\item[(iv)]  $\adfPENT(4,r)$ for $r \in $
\{13, 17, 20, 21, 24,
29, 33, 37, 40, 45, 49, 52, 53, 60, 61, 65, 69, 77, 80, 81, 85, 93, 97, 100, 101, 108, 109, 117, 120, 125, 133, 140, 141, 149, 157, 160, 165, 173, 180\}, \cite{Forbes2020P};
\item[(v)]  $\adfPENT(5,r)$ for $r \in $ \{20, 25, 30, 35, 40\}, \cite{Forbes2020P};
\item[(vi)] the $\adfPENT(6,7)$ and $\adfPENT(7,7)$ based on the Hoffman--Singleton graph, \cite{BallBambergDevillersStokes2013}.
\end{enumerate}

In the next section we describe the basic properties of generalized pentagonal geometries,
relating them to the corresponding results of \cite{BallBambergDevillersStokes2013}.
We find that although many of the lemmas in \cite[Section 2]{BallBambergDevillersStokes2013} concerning pentagonal geometries
transfer straightforwardly to $\adfPENT(k,r,w)$, there are some properties that do not.
For example, when $r > (w - 1)/(k - 1)$ a $\adfPENT(k,r,w)$ can have a connected deficiency graph with girth 4
whereas this is impossible if $w = k$.
We also settle the existence problem for $\adfPENT(2,r,w)$.
In Section 3 we show how to construct $\adfPENT(3, r, w)$ geometries with deficiency $(w,5^+)$-graphs,
and we prove that these exist for $w \in$ \{7, 9, 13, 15, 19, 21, 25, 27, 31, 33\} and all sufficiently large admissible $r$.
In Section 4 we present a number of $\adfPENT(4, r, 13)$ geometries with connected deficiency $(13,5^+)$-graphs.
Consequently we show that a $\adfPENT(4, r, 13)$ with deficiency $(13,5^+)$-graph exists for every sufficiently large even $r$.


\section{Generalized pentagonal geometries}\label{sec:Generalized pentagonal geometries}

\begin{figure}[h!]
\begin{center}
\includegraphics[width=80mm,trim=0mm 0mm 0mm 0mm, clip]{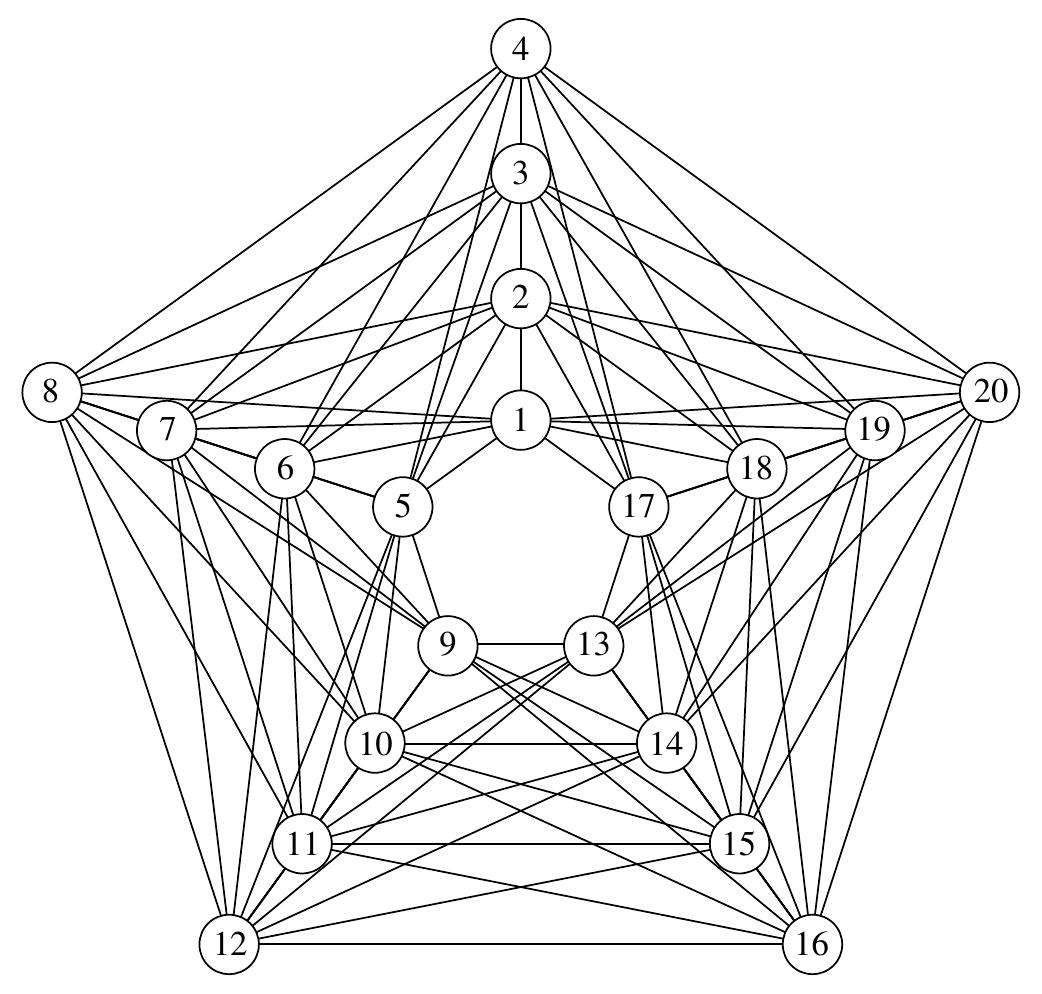} 
\end{center}
\caption{$\adfPENT(2, 3m - 1, 2m)$, $m = 4$}\label{fig:PENT-2-r-w}
\end{figure}

%
\begin{lemma} \label{lem:v-b}
A generalized pentagonal geometry $\adfPENT(k, r, w)$ has $r(k-1) + w + 1$ points and $(r(k - 1) + w + 1)r/k$ lines.
A $\adfPENT(k, r, w)$ exists only if $r(r - w - 1) \equiv 0 \adfmod{k}$.
\end{lemma}
\begin{proof}
Consider a point, $x$. Every other point is either on one of the $r$ lines that contain $x$ or is one of the $w$ points of $x^{\mathrm{opp}}$;
$r(k-1) + w + 1$ points altogether.

Let $v = r(k-1) + w + 1$. The number of `point--lines' is $vr$. However, each line contains $k$ points.
Therefore the number of lines is $vr/k$, which must of course be an integer.
\end{proof}
Given $k$ and $w$, we say that $r$ is {\em admissible} if $r(w + 1 - r) \equiv 0 \adfmod{k}$.
%
\begin{lemma} \label{lem:opp-des-disjoint}
If $x$ and $y$ are distinct points in a generalized pentagonal geometry such that there exist two distinct points common to both
$x^\mathrm{opp}$ and $y^\mathrm{opp}$, then $x^\mathrm{opp}$ and $y^\mathrm{opp}$ have a common block.
\end{lemma}
\begin{proof}
If $a$ and $b$ are points in both $x^\mathrm{opp}$ and $y^\mathrm{opp}$,
then the block containing $a$ and $b$ must belong to both opposite designs.
\end{proof}
In \cite[Lemma 2.3]{BallBambergDevillersStokes2013} it is shown that if $x$ and $y$ are distinct points of a
$\adfPENT(k,r)$ such that $x^\mathrm{opp} = y^\mathrm{opp} = L$,
then there exists a unique opposite line pair $\{L, M\}$ such that $x, y \in M$.
However, this does not generally extend to geometries $\adfPENT(k,r,w)$ where $w > k$.
For example, take the $\adfPENT(2,11,8)$ shown in Figure~\ref{fig:PENT-2-r-w}, which is clearly a
special case of a general construct, $\adfPENT(2, 3m - 1, 2m)$.
Here we have five copies of $K_4$, represented as the `spokes', linked by the edges of five copies of $K_{4,4}$ and, for instance,
$1^\mathrm{opp} = 2^\mathrm{opp} = 3^\mathrm{opp} = 4^\mathrm{opp} = \{9, 10, \dots, 16\}$ whereas
$x^\mathrm{opp} \neq y^\mathrm{opp}$ for $x \in \{9, 10, 11, 12\}$, $y \in \{13, 14, 15, 16\}$.
Observe also that the deficiency graph is connected, is 8-regular and has girth 4.
%
\begin{lemma} \label{lem:deficiency-graph}
Let $D$ be the deficiency graph of a generalized pentagonal geometry $\adfPENT(k,r,w)$.
Then $D$ is $w$-regular and has girth at least $4$.

For distinct points $x$ and $y$, let $U_{x,y}$ denote the set of points common to both $x^\mathrm{opp}$ and $y^\mathrm{opp}$,
and let $u_{x,y} = |U_{x,y}|$.
Then $U_{x,y}$ is the point set of a Steiner system $S(2,k,u_{x,y})$ that exists as a
subdesign of both $x^\mathrm{opp}$ and $y^\mathrm{opp}$.
\end{lemma}
\begin{proof}
For any vertex $x$ of $D$, the neighbourhood of $x$ is the point set of $x^\mathrm{opp}$, which has cardinality $w$.
If $D$ contains a path, $a \sim x \sim b$, then $a$ and $b$ are in $x^\mathrm{opp}$ and must therefore be collinear.
Hence $D$ is triangle-free.
The neighbours of vertices $x$ and $y$ in $D$ are clearly the point sets of $x^\mathrm{opp}$ and $y^\mathrm{opp}$ respectively.

When $u_{x,y} = 0$ it follows vacuously that an $S(2,k,0)$ is a subdesign of both $x^\mathrm{opp}$ and $y^\mathrm{opp}$.
If $u_{x,y} = 1$, then $U_{x,y}$ is clearly the point set of an $S(2,k,1)$ which is a subdesign of both $x^\mathrm{opp}$ and $y^\mathrm{opp}$.

Suppose $u_{x,y} \ge 2$ and note that the points of $U_{x,y}$ are pairwise collinear.
Hence, for any pair of distinct points $a$ and $b$, there is a line containing both $a$ and $b$ that must belong to
the line set of both $x^\mathrm{opp}$ and $y^\mathrm{opp}$.
Let $\mathcal{L}$ be the set of all these lines.
Then $(U_{x,y}, \mathcal{L})$ is a Steiner system $S(2,k,u_{x,y})$, and
$\mathcal{L}$ is contained in the line sets of both $x^\mathrm{opp}$ and $y^\mathrm{opp}$.
\end{proof}
The possible values of $u_{x,y}$ in Lemma~\ref{lem:deficiency-graph} always include 0, 1, $k$ and $w$,
and by Fisher's inequality for block designs, if $w > k$, then $w \ge k^2 - k +1$.
Clearly we cannot have $1 < u_{x,y} < k$ whereas
for values strictly between $k$ and $w$, it depends on whether there exists a Steiner system $S(2,k,u_{x,y})$
imbedded in both of the Steiner systems $S(2,k,w)$ that form $x^\mathrm{opp}$ and $y^\mathrm{opp}$.
For example, in a $\adfPENT(3,r,19)$ where the STS(19)s used for the opposite designs include subdesigns STS(7) and STS(9),
\cite[Table 4]{ColbournForbesGrannellGriggsKaskiOstergardPikePottonen2010},
one can have $u_{x,y} \in \{0, 1, 3, 7, 9, 19\}$.
%
\begin{lemma} \label{lem:pairwise disjoint opposite designs}
For a generalized pentagonal geometry $\adfPENT(k, r, w)$, the following properties are equivalent:
\begin{enumerate}
\item[(i)] for any two distinct points $x$ and $y$, $x^\mathrm{opp}$ and $y^\mathrm{opp}$ have no common block;
\item[(ii)] for any two distinct points $x$ and $y$, the point sets of $x^\mathrm{opp}$ and $y^\mathrm{opp}$ have at most one common point;
\item[(iii)] the deficiency graph has girth at least $5$.
\end{enumerate}
\end{lemma}
\begin{proof}
(i) $\Rightarrow$ (ii): Suppose for distinct points $x$ and $y$, $x^\mathrm{opp}$ and $y^\mathrm{opp}$ have distinct common points $a$ and $b$.
Then there exists a unique line containing $a$ and $b$ that must be present in the line sets of both $x^\mathrm{opp}$ and $y^\mathrm{opp}$.

(ii) $\Rightarrow$ (iii): We know from Lemma~\ref{lem:deficiency-graph} that the deficiency graph has girth at least 4.
Suppose the deficiency graph contains a 4-cycle, $(x, a, y, b)$.
Then $a$ and $b$ are in both $x^\mathrm{opp}$ and $y^\mathrm{opp}$.

(iii) $\Rightarrow$ (i): Suppose for distinct points $x$ and $y$, $x^\mathrm{opp}$ and $y^\mathrm{opp}$
have a common block containing points $a$ and $b$.
Then the deficiency graph contains a 4-cycle, $(x, a, y, b)$.
\end{proof}
For clarification in the rest of the paper,
we define the property {\em pairwise disjoint opposite designs} of a generalized pentagonal geometry
to mean that for any two distinct points $x$ and $y$, $x^\mathrm{opp}$ and $y^\mathrm{opp}$ have no common block,
i.e.\ (i) of Lemma~\ref{lem:pairwise disjoint opposite designs}.
Thus we are referring to the whole collection of opposite designs, including possible repetitions, and
we allow two opposite designs to have at most one common point.
As we have seen, the property is equivalent to having a deficiency graph with girth at least 5.
%
\begin{lemma} \label{lem:r >= w(w-1)/(k-1)}
For a $\adfPENT(k, r, w)$ with deficiency graph of girth at least $5$, we have $r \ge w(w - 1)/(k - 1)$.
\end{lemma}
\begin{proof}
Let $v$ be the number of points in the $\adfPENT(k, r, w)$.
An opposite design contains $w(w - 1)/(k(k - 1))$ blocks.
Since the opposite designs are pairwise disjoint by Lemma~\ref{lem:pairwise disjoint opposite designs},
the number of blocks, $vr/k$, cannot be less than $vw(w - 1)/(k(k - 1))$.
\end{proof}
When equality occurs in Lemma~\ref{lem:r >= w(w-1)/(k-1)} the number of points is $v = w^2 + 1$,
the Moore bound for $(w,5)$-graphs.
%
\begin{lemma} \label{lem:deficiency-graph-distance-2}
Let $x$ and $y$ be two distinct points in a generalized pentagonal geometry with deficiency graph $D$.
Then the distance in $D$ between $x$ and $y$ is $2$ if and only if there exists a point $z$ such that
$z^\mathrm{opp}$ contains both $x$ and $y$.
\end{lemma}
\begin{proof}
This is similar to \cite[Lemma 2.7]{BallBambergDevillersStokes2013}. Suppose $x$ and $y$ are distinct points.
If $x$ and $y$ are at distance 2 in $D$, then $x \not\sim y$ and there exists a point $z$ such that $x \sim z \sim y$;
hence $x$ and $y$ belong to $z^\mathrm{opp}$.
On the other hand, if $x$ and $y$ belong to $z^\mathrm{opp}$, then there is a line in $z^\mathrm{opp}$ containing both $x$ and $y$;
therefore $x \sim z \sim y$ and $x \not\sim y$.
\end{proof}
Thus the deficiency graph of a $\adfPENT(k, r, w)$ with pairwise disjoint opposite designs
has $r(k - 1) + w + 1$ vertices and is a $(w,5^+)$-graph.
Also the distance between vertices $x$ and $y$ is at least 3 if $x$ and $y$ are distinct points
on the same line that is not part of an opposite design.

%
\begin{theorem} \label{thm:PENT(2,r,w)-exist}
Suppose $r \ge w - 1 \ge 1$ and $rw$ is even. Let $v = r + w + 1$.
A generalized pentagonal geometry $\adfPENT(2,r,w)$ exists
if and only if there exists a $w$-regular graph with $v$ vertices and girth at least $4$.
A generalized pentagonal geometry $\adfPENT(2,r,w)$ with pairwise disjoint opposite designs exists
if and only if there exists a $w$-regular graph with $v$ vertices and girth at least $5$.
\end{theorem}
\begin{proof}
The inequality for $r$ ensures that $v \ge 2w$, the Moore bound for a $(w,4)$-graph.
Suppose $D$ is a $v$-vertex $(w,4^+)$-graph.
Then $wv \equiv rv \equiv 0 \adfmod{2}$
and the edges of the complement of $D$ are clearly the $v(v-1)/2 - wv/2 = vr/2$ blocks of a $\adfPENT(2,r,w)$
in which for each point $x$, the neighbours of vertex $x$ in $D$ comprise the point set of the Steiner system $S(2,2,w)$
that forms $x^\mathrm{opp}$.

Furthermore, if there exist points $x$ and $y$ such that $x^\mathrm{opp}$ and $y^\mathrm{opp}$ both contain a block $\{a,b\}$, say,
then there is a 4-cycle, $(x,a,y,b)$.
Hence the opposite designs must be pairwise disjoint if $D$ has girth at least 5.

Conversely, by Lemmas~\ref{lem:deficiency-graph} and \ref{lem:pairwise disjoint opposite designs}, the deficiency graph of a $\adfPENT(2,r,w)$ has the stated properties.
\end{proof}
If $w = k$ and $x^\mathrm{opp} = y^\mathrm{opp}$ for distinct points $x$ and $y$,
then $x^\mathrm{opp}$ and the line containing $\{x,y\}$ must form an opposite line pair, \cite[Lemma 2.3]{BallBambergDevillersStokes2013}.
On the other hand, when $w > k$ the situation is not so simple.
As we have seen, it is possible to have distinct points $x$, $y$ such that $x^\mathrm{opp}$ and $y^\mathrm{opp}$ contain a common block
but do not form an opposite design pair.
We leave the investigation of this complication as well as a detailed discussion of geometries containing opposite design pairs
for a possible future paper.
Henceforth we shall confine our attention mainly to generalized pentagonal geometries where the opposite designs are pairwise disjoint,
or equivalently by Lemma~\ref{lem:pairwise disjoint opposite designs},
where the deficiency graph has girth at least 5.
Now pairs of vertices of the deficiency graph of a $\adfPENT(k,r,w)$ with $v = (k - 1)r + w + 1$
may be categorized according to their contribution to the geometry:
\begin{enumerate}
\item[(i)] $v w/2$ adjacent pairs, which necessarily have no common neighbour, are not collinear points in the geometry;
\item[(ii)] $w(w - 1)v/2$ non-adjacent pairs with one common neighbour are collinear in the blocks of the opposite designs;
\item[(iii)] $(v - w^2 - 1)v/2$ non-adjacent pairs with no common neighbour are collinear in the blocks that do not occur in opposite designs.
\end{enumerate}
Another way of looking at a $\adfPENT(k,r,w)$ is given by the next theorem.
%
\begin{theorem}
\label{thm:edgewise-decomposition}
Let $k \ge 2$, $w \ge k$ and $r \ge w(w - 1)/(k - 1)$ be integers.
Let $v = (k-1)r + w + 1$ and $V = \{1, 2, \dots, v\}$.
Then there exists a pentagonal geometry $\adfPENT(k,r,w)$ with point set $V$ and having a deficiency graph of girth at least $5$
if and only if there exists an edgewise decomposition of the complete graph $K_v$ with vertex set $V$ into
the following graphs, each with vertex set $V$:
\begin{enumerate}
\item[(i)] a $w$-regular graph, $D$, with girth at least $5$;
\item[(ii)] a $w(w - 1)$-regular graph $S$ which itself has an edgewise decomposition into
$v$ complete graphs $K_w$, $S_1$, $S_2$, \dots, $S_v$, say, such that for $i = 1$, $2$, \dots, $v$,
the vertex set of $S_i$ is the neighbourhood set of vertex $i$ of $D$, and
$S_i$ admits an edgewise decomposition into $w(w - 1)/(k(k - 1))$ complete graphs $K_k$;
\item[(iii)] if $r > w(w - 1)/(k - 1)$, a $(v - w^2 - 1)$-regular graph, $T$,
which admits an edgewise decomposition into $(v - w^2 - 1)v/(k(k - 1))$ complete graphs $K_k$.
Moreover, if $r = w(w - 1)/(k - 1) + 1$, the graphs $K_k$ partition the vertex set of $T$.
\end{enumerate}
\end{theorem}
\begin{proof}
A Steiner system $S(2,k,w)$ exists if and only if
there exists an edgewise decomposition of the complete graph $K_w$ into complete graphs $K_k$.

The neighbourhood set of vertex $x$ of $D$ is the point set of $x^\mathrm{opp}$, and
$x^\mathrm{opp}$ is the $S(2,k,w)$ corresponding to $S_x$,

Since $D$ is a $(w, 5^+)$ graph, two neighbourhood sets of $D$ cannot have more than one common vertex, and
a given vertex $x$ of $D$ occurs in precisely $w$ of the graphs $S_i$, i.e.\ those corresponding to $\{y^\mathrm{opp}: y \in x^\mathrm{opp}\}$.
Moreover, considered as a point in the Steiner system $S(2,w,k)$, $x$ occurs in precisely $(w - 1)/(k - 1)$ of its blocks.
Therefore $x$ occurs in precisely $w - 1$ edges of those $S_i$ that contain $x$.
Hence $S$ is $w(w - 1)$-regular and $T$ is $(v - 1 - w - w(w - 1))$-regular.
Also $T$ has an edgewise decomposition into $K_k$ graphs corresponding to the
$$\dfrac{v}{k} \left(r - \dfrac{w(w - 1)}{k - 1}\right) = \dfrac{v}{k} \cdot \dfrac{v - w^2 - 1}{k - 1}$$
lines that do not belong to any opposite design.
If $r = w(w - 1)/(k - 1) + 1$, then $T$ is $(k - 1)$-regular and
has an edgewise decomposition into $v/k$ $K_k$ graphs, which must therefore partition $V$.
\end{proof}
%
\begin{theorem} \label{thm:minimum-r}
$\mathrm{(i)}$ There is no $\adfPENT(k, w(w - 1)/(k - 1), w)$ with deficiency graph of girth $5$ if
$w \ge k \ge 2$ and $w \notin \{2, 3, 7, 57\}$.

$\mathrm{(ii)}$ Suppose $w \ge k \ge 2$, $w \in \{2, 3, 7\}$ and $r = w(w - 1)/(k - 1)$.
Then there exists a $\adfPENT(k, r, w)$ with deficiency graph of girth $5$ if and only if
\begin{equation}\label{eqn:(k-r-w)}
(k,r,w) \in \{(2,2,2), (2,6,3), (3,3,3), (2,42,7), (3,21,7), (7,7,7)\}.
\end{equation}
\end{theorem}
\begin{proof}
(i) The number of points is $w^2 + 1$, the Moore bound for $(w,5)$-graphs.
Hence a $\adfPENT(k, w(w - 1)/(k - 1), w)$ cannot exist unless $w$ has one of the stated values;
see, for example, \cite[Section 11.12]{Cameron1994}.

(ii) The existence of a $\adfPENT(k, w(w - 1)/(k - 1), w)$ requires a Steiner system $S(2,k,w)$,
and therefore the only possibilities are given by (\ref{eqn:(k-r-w)}).
For $(2,2,2)$, $(3,3,3)$ and $(7,7,7$), see \cite[Section 3]{BallBambergDevillersStokes2013}.
For $(2,6,3)$ take the Petersen graph and overlay the neighbourhood set of each point with a triangle.
For $(2,42,7)$ and $(3,21,7)$, take the Hoffman--Singleton graph and overlay the neighbourhood set of each point with a $K_7$
or an STS(7), as appropriate.
\end{proof}
If the 3250-vertex, 57-regular Moore graph with girth 5 exists,
then so do generalized pentagonal geometries which have it as the deficiency graph:
$\adfPENT(57, 57, 57)$, \cite{BallBambergDevillersStokes2013},
as well as $\adfPENT(2,3192,57)$, $\adfPENT(3,1596,57)$ and $\adfPENT(8,456,57)$
obtained by overlaying the neighbourhood sets of the graph with
complete graphs $K_{57}$, STS(57)s and projective planes of order 7 respectively.
%
\begin{theorem} \label{thm:minimum-r-plus-1}
There is no $\adfPENT(k, w(w - 1)/(k - 1) + 1, w)$ with deficiency graph of girth at least $5$
for $w \ge k \ge 2$ except for $\adfPENT(2, 3, 2)$, $\adfPENT(6, 7, 6)$ and possibly $\adfPENT(56, 57, 56)$.
\end{theorem}
\begin{proof}
For $w = k$, see \cite[Corollary 4.5]{BallBambergDevillersStokes2013}.

Let $w > k \ge 2$, $r = w(w - 1)/(k - 1) + 1$ and
suppose there exists a $\adfPENT(k, r, w)$ with deficiency $(w,5^+)$-graph.
The number of points and the number of lines are respectively
$$v = w^2 + k \text{\quad and \quad} b = \dfrac{v r}{k} = \dfrac{(w^2 + k) (w^2 - w + k - 1)}{k(k - 1)}.$$
Since the Moore bound for $(w,6)$-graphs, namely $2(w^2 - w + 1)$, exceeds $w^2 + k$ when $w > k \ge 2$,
we may assume that the deficiency graph has girth 5.

In \cite{Brown1967} it is proved that the smallest number of vertices for which a $(w,5)$-graph exists is not equal to $w^2 + 2$.
Since $w > 2$, the Moore bound, $w^2 + 1$, is not attained when $w$ is even and
consequently there is no $(w,5)$-graph with $w^2 + 2$ vertices.
Hence there is no $\adfPENT(2, r, w)$ with $w > 2$.
In what follows we assume $k \ge 3$.

For $k = 3$, it is known that there is no $(w,5)$-graph with $w^2 + 3$ vertices except possibly when
$w$ is even or when $w$ is odd and $w = \ell^2 + \ell + 1 \pm 2$ for some integer $\ell$, \cite{Kovacs1981}.
Although it would suffice to deal only with these exceptional $w$, it is convenient to include $k = 3$ in our proof.

The number of lines in opposite designs is $v w(w - 1)/(k(k-1))$ and therefore
the number of non-opposite lines, i.e.\ lines that do not occur in opposite designs, is
$$b_\mathrm{rem} = b - \dfrac{v w(w - 1)}{k(k-1)} = \dfrac{v}{k}.$$
Hence $k \mid v$ and $k \mid w^2$.
Moreover, if a Steiner system $S(2, k, w)$ exists, then $k(k - 1) \mid w(w - 1)$ and $k - 1 \mid w - 1$.
It follows that $k \mid w$ and
\begin{equation}\label{eqn:minimum-w-t-k}
w = t k (k - 1) + k,~~~~~ t = 1, 2, \dots.
\end{equation}
Let $D$ denote the deficiency graph and let $E$ denote the distance at least 3 graph,
where the vertices of $E$ are points of the geometry and there is an edge $i \sim j$
whenever the distance between $i$ and $j$ in $D$ is at least 3.
Let $A$ and $B$ be the adjacency matrices of $D$ and $E$ respectively and
denote the distance in $D$ between vertices $i$ and $j$ by $d(i,j)$.
Then, recalling that $[A^2]_{i,j}$ is the number of 2-step walks from vertex $i$ to vertex $j$,
in row $i$ of $A^2$ we have
\begin{equation}\label{eqn:minimum-r-plus-1-freq-A}
\left\{\begin{array}{lll}
\left[A^2\right]_{i,i}                &= w,&~~~ \textrm{which~occurs~once},\\
\left[A^2\right]_{i,j,~i \sim j}      &= 0,&~~~ \textrm{which~occurs $w$ times},\\
\left[A^2\right]_{i,j,~ d(i,j) = 2}   &= 1,&~~~ \textrm{which~occurs $w(w - 1)$ times},\\
\left[A^2\right]_{i,j,~ d(i,j) \ge 3} &= 0,&~~~ \textrm{which~occurs $k-1$ times}.
\end{array} \right.
\end{equation}
Let $I$, $J$ and $\bold{j}$ denote the identity matrix, the all-ones square matrix and the all-ones vector respectively, each of dimension $v$.
From (\ref{eqn:minimum-r-plus-1-freq-A}) it is plain that $E$ is $(k - 1)$-regular and
\begin{equation}\label{eqn:minimum-r-plus-1-B}
B = J - A^2 - A + (w - 1)I.
\end{equation}
By Lemma~\ref{lem:deficiency-graph-distance-2}, two distinct vertices of $E$ are
adjacent precisely when they are on the same non-opposite line.
Hence $E$ consists of $b_\mathrm{rem} = v/k$ pairwise disjoint complete graphs $K_k$.
Since $D$ is $w$-regular and is too small to have more than one component,
$A$ has eigenvalue $w$ with multiplicity 1 and corresponding eigenvector $\bold{j}$.
Let $\bold{s}$ be an eigenvector of $A$ that is orthogonal to $\bold{j}$ and let $s$ be its eigenvalue.
Then from (\ref{eqn:minimum-r-plus-1-B}) we have
$$B\,\bold{s} = (J - A^2 - A + (w - 1)I)\,\bold{s} =  (- s^2 - s + w - 1)\,\bold{s},$$
and hence $- s^2 - s + w - 1$ is an eigenvalue of $B$.
But the eigenvalues of $B$ are $k - 1$ and $-1$.
Therefore the eigenvalues of $A$ comprise:
\begin{enumerate}
\item[(i)] $w$ with multiplicity 1 and eigenvector $\bold{j}$;
\item[(ii)] roots of $-s^2 - s + w - 1 = k - 1$,
       $$s_1 = \dfrac{-1 + \sqrt{4 w - 4k + 1}}{2},~~ s_2 = \dfrac{-1 - \sqrt{4 w - 4k + 1}}{2},$$
       with multiplicities $m_1$, $m_2$ respectively;
\item[(iii)] roots of $-s^2 - s + w - 1 = -1$,
       $$s_3 = \dfrac{-1 + \sqrt{4 w + 1}}{2},~~ s_4 = \dfrac{-1 - \sqrt{4 w + 1}}{2},$$
       with multiplicities $m_3$, $m_4$ respectively.
\end{enumerate}
Since $A$ is the adjacency matrix of a triangle-free graph, both $A$ and $A^3$ have zero trace.
Moreover, by (\ref{eqn:minimum-r-plus-1-freq-A}), the trace of $A^2$ is $v w$.
We therefore have
$$1 + \sum_{i=1}^4 m_i = v,~~~ w^2 + \sum_{i=1}^4 m_i s_i^2 = v w,~~~ v = w^2 + k,$$
$$w + \sum_{i=1}^4 m_i s_i = w^3 + \sum_{i=1}^4 m_i s_i^3 = 0,$$
which has the solution
\begin{align*}
m_1 &= \dfrac{w(w - 2 k  + w \sqrt{1 - 4 k + 4 w})}{2 k \sqrt{1 - 4 k + 4 w}},\\
m_2 &= \dfrac{w(2 k - w  + w \sqrt{1 - 4 k + 4 w})}{2 k \sqrt{1 - 4 k + 4 w}},\\
m_3 &= \dfrac{(\sqrt{4 w + 1} + 1) (k^2 + (w^2 - 1) k - w^2)}{2 k \sqrt{4 w + 1}},\\
m_4 &= \dfrac{(\sqrt{4 w + 1} - 1) (k^2 + (w^2 - 1) k - w^2)}{2 k \sqrt{4 w + 1}}.
\end{align*}
By (\ref{eqn:minimum-w-t-k}), we cannot have $w = 2k$ or $k^2 + (w^2 - 1) k - w^2) = 0$;
so if one of the eigenvalues is irrational, then at least one of the multiplicities is irrational.

Hence we may assume all eigenvalues of $A$ are rational and there exist integers $p$ and $q$ such that
$$\sqrt{4 w - 4k + 1} = 2p + 1,~~~ \sqrt{4 w + 1} = 2q + 1.$$
Therefore
$$w = p^2 + p + k = q^2 + q$$
and it follows that $(q - p)(q + p + 1) = k = de$, say, where $d = q - p$ is a divisor of $k$.
Then
$$p = \dfrac{e - d - 1}{2},~~~ q = \dfrac{e + d - 1}{2}$$
and
$$w - k = p^2 + p = \dfrac{(d - e)^2 - 1}{4} < tk(k - 1) = w - k,$$
a contradiction.
\end{proof}

For creating new generalized pentagonal geometries from existing ones, we use Theorem~\ref{thm:GDD-basic}, below,
a straightforward adaptation of Wilson's Fundamental Construction, \cite{WilsonRM1972}, \cite[Theorem IV.2.5]{GreigMullen2007}.
For this theorem, the following definition is central.

A {\em group divisible design}, $k$-GDD, of type $g_1^{u_1} g_2^{u_2} \dots g_n^{u_n}$ is
an ordered triple ($V,\mathcal{G},\mathcal{B}$) such that:
\begin{adfenumerate}
\item[(i)]{$V$
is a base set of cardinality $u_1 g_1 + u_2 g_2 + \dots + u_n g_n$;}
\item[(ii)]{$\mathcal{G}$
is a partition of $V$ into $u_i$ subsets of cardinality $g_i$, $i = 1, 2, \dots, n$, called \textit{groups};}
\item[(iii)]{$\mathcal{B}$
is a non-empty collection of $k$-subsets of $V,$ called \textit{blocks}; and}
\item[(iv)]{each pair of elements from distinct groups occurs in precisely one block but no pair of
elements from the same group occurs in any block.}
\end{adfenumerate}
A Steiner system $S(2,k,w)$ is essentially a $k$-GDD of type $1^w$.
A $k$-GDD of type $q^k$ is also called a {\em transversal design}, TD$(k,q)$.

%
\begin{theorem}
\label{thm:GDD-basic}
Let $k \ge 2$ and $w \ge k$ be integers.
For $i = 1$, $2$, \dots, $n$, let $r_i$ be a positive integer, let $v_i = (k-1) r_i + w + 1$,
and suppose there exists a generalized pentagonal geometry $\adfPENT(k, r_i, w)$.
Suppose also that there exists a $k$-$\mathrm{GDD}$ of type $v_1^{u_1} v_2^{u_2} \dots v_n^{u_n}$.
Let $N = u_1 + u_2 + \dots + u_n$ and
$R = u_1 r_1 + u_2 r_2 + \dots + u_n r_n$.
Then there exists a generalized pentagonal geometry
$\adfPENT(k, R + (N - 1)(w + 1)/(k - 1), w)$.

Furthermore, if for $i = 1$, $2$, \dots, $n$, the $\adfPENT(k, r_i, w)$ has a deficiency graph of girth at least $5$,
then so does the $\adfPENT(k, R + (N - 1)(w + 1)/(k - 1), w)$.
\end{theorem}
\begin{proof}
Overlay each group of size $v_i$ with a $\adfPENT(k, r_i, w)$, $i = 1$, 2, \dots, $n$.
The total number of points in the geometry that results is $v_1 u_1 + v_2 u_2 + \dots + v_n u_n = R(k - 1) + N(w + 1)$.
If the $\adfPENT(k, r_i, w)$, $i = 1$, 2, \dots, $n$, have deficiency $(w,5^+)$-graphs,
then by Lemma~\ref{lem:pairwise disjoint opposite designs} they have pairwise disjoint opposite designs,
which are clearly not affected by the blocks of the group divisible design.
Therefore, again by Lemma~\ref{lem:pairwise disjoint opposite designs}, the deficiency graph of the constructed geometry has girth at least 5.
\end{proof}
%


\section{Block size 3}\label{sec:Block size 3}
By Lemma~\ref{lem:v-b}, a pentagonal geometry $\adfPENT(3,r,w)$ has $v = 2r + w + 1$ points and $r(2r + w + 1)/3$ lines.
Hence $r$ is admissible if and only if $r \equiv 0 \textrm{~or~} w + 1 \adfmod{3}$.
A Steiner triple system of order $w \ge 0$ exists if and only if $w \equiv 1 \textrm{~or~} 3 \adfmod{6}$.
Furthermore, $v$ is even,
and for constructions using Theorem~\ref{thm:GDD-basic} we employ the following
lemma concerning the existence of group divisible designs with block size 3 and even group sizes.
%
\begin{lemma}
\label{lem:3-GDD-existence}
Suppose $g \ge 1$ and $g$ is even. Then there exists a {\rm 3-GDD} of type $g^u$ if $u \ge 3$ and $u(u-1)g \equiv 0 \adfmod{3}$.

Suppose also that $m$ is even and $m \le g(u-1)$.
Then there exists a {\rm 3-GDD} of type $g^u m^1$ if $u \ge 3$ and $g^2 u (u-1) + 2 g u m \equiv 0 \adfmod{3}$.
\end{lemma}
\begin{proof}
See \cite{ColbournHoffmanRees1992} or \cite[Theorem IV.4.2]{Ge2007}.
\end{proof}
%
\begin{lemma}
\label{lem:PENT-3-r-7-direct}
There exists a $\adfPENT(3,r,7)$ with connected deficiency graph of girth at least $5$ for
$33 \le r \le 149$, $r \equiv 0 \textrm{~or~} 2 \adfmod{3}$.
\end{lemma}
\begin{proof}
Since these geometries are not difficult to find we shall, as a reasonable alternative to presenting each one in detail,
merely describe how they are obtained. Let $v = 2r + 8$.

By some means create a random connected $(7,5^+)$-graph with $v$ vertices.
The points of the geometry are the vertices of the graph.
For each vertex $x$, construct a Steiner triple system of order 7 on the neighbours of $x$.
These form the blocks of $x^\mathrm{opp}$.
Finally, create the remaining blocks of the geometry by hill climbing, as described in \cite{Stinson1985} or \cite[Section 2.7.2]{ColbourRosa1999}.
Of course this last step is doomed if the chosen graph is not the deficiency graph of a generalized pentagonal geometry.
We must therefore be prepared to restart from the beginning with a different graph if it becomes evident that the hill climbing
process does not appear to be working.

We give two examples.
In each case we specify $d$ followed by $dr$ numbers.
To obtain the geometry, gather the numbers into triples to form a set of $dr/3$ base blocks, $B$,
which are then developed into the line set of the geometry by the mapping $x \mapsto x + d \adfmod{v}$.
The point set is  $\{0, 1, \dots, v - 1\}$ and
the first $d$ sets of 7 triples of $B$ are $0^\mathrm{opp}$, $1^\mathrm{opp}$, \dots, $(d-1)^\mathrm{opp}$.
Further examples are given in the Appendix.


{\boldmath $\adfPENT(3,33,7)$}, $d = 2$:
{5, 24, 49;  5, 50, 51;  5, 59, 63;  24, 50, 59;  24, 51, 63;  49, 50, 63;  49, 51, 59;  12, 16, 19;  12, 24, 26;  12, 57, 70;  16, 24, 57;  16, 26, 70;  19, 24, 70;  19, 26, 57;  0, 6, 42;  0, 11, 53;  0, 15, 65;  0, 17, 57;  0, 18, 52;  0, 19, 71;  0, 21, 47;  0, 37, 43}

{\boldmath $\adfPENT(3,35,7)$}, $d = 6$:
{11, 19, 22;  11, 58, 59;  11, 65, 74;  19, 58, 65;  19, 59, 74;  22, 58, 74;  22, 59, 65;  2, 31, 33;  2, 35, 39;  2, 49, 60;  31, 35, 49;  31, 39, 60;  33, 35, 60;  33, 39, 49;  1, 6, 21;  1, 29, 45;  1, 46, 52;  6, 29, 46;  6, 45, 52;  21, 29, 52;  21, 45, 46;  38, 40, 43;  38, 49, 58;  38, 62, 70;  40, 49, 62;  40, 58, 70;  43, 49, 70;  43, 58, 62;  15, 24, 27;  15, 32, 38;  15, 45, 60;  24, 32, 45;  24, 38, 60;  27, 32, 60;  27, 38, 45;  18, 23, 24;  18, 49, 56;  18, 65, 72;  23, 49, 65;  23, 56, 72;  24, 49, 72;  24, 56, 65;  26, 0, 61;  18, 31, 19;  44, 74, 32;  60, 42, 46;  19, 45, 77;  17, 20, 69;  28, 49, 4;  45, 64, 3;  44, 7, 61;  76, 11, 38;  43, 53, 39;  8, 18, 61;  52, 24, 36;  33, 32, 11;  24, 31, 16;  18, 27, 70;  26, 29, 8;  68, 24, 58;  59, 64, 69;  29, 9, 58;  24, 51, 26;  63, 13, 49;  4, 6, 23;  31, 9, 53;  0, 10, 75;  0, 35, 37;  2, 16, 41;  3, 41, 53}
\end{proof}
%
\begin{theorem}
\label{thm:PENT-3-r-7-constructed}
There exists a $\adfPENT(3,r,7)$ with deficiency graph of girth at least $5$ if and only if $r \ge 21$ and $r \equiv 0 \textrm{~or~} 2 \adfmod{3}$,
except possibly for $r \in $
\{$23$, $24$, $26$, $27$, $29$, $30$, $32$\}.
\end{theorem}
\begin{proof}

The smallest $(7,5)$-graph is the 50-vertex Hoffman--Singleton graph,
which yields a $\adfPENT(3,21,7)$ by Theorem~\ref{thm:minimum-r}.
Therefore, by Lemma~\ref{lem:PENT-3-r-7-direct}, we need to consider only admissible $r \ge 150$.

Let $S_7 = \{s: 33 \le s \le 71,\; s \equiv 0 \textrm{~or~} 2 \adfmod{3}\}$ and observe that
the elements of $S_7$ cover the admissible residues modulo 39.
By Lemma~\ref{lem:PENT-3-r-7-direct}, there exists a $\adfPENT(3,s,7)$ for $s \in S_7$.

For each $s \in S_7$, take $u$ copies of a $\adfPENT(3,35,7)$, which has 78 points,
and one copy of a $\adfPENT(3,s,7)$, which has $2s + 8$ points.
By Lemma~\ref{lem:3-GDD-existence}, since $2 s + 8 \le 78\cdot 2$,
there exists a 3-GDD of type $78^u (2s + 8)^1$ for $u = 3$, 4, \dots.
Now use Theorem~\ref{thm:GDD-basic} to construct a $\adfPENT(3,39u + s,7)$ for $u \ge 3$.

Hence there exists a $\adfPENT(3,r,7)$ for admissible $r \ge 150$.
\end{proof}
%
\begin{lemma}
\label{lem:PENT-3-r-9-direct}
There exists a $\adfPENT(3,r,9)$ with connected deficiency graph of girth at least $5$ for
$57 \le r \le 244$, $r \equiv 0 \textrm{~or~} 1 \adfmod{3}$.
\end{lemma}
\begin{proof}

We give two examples, below,
in the same format as for Lemma~\ref{lem:PENT-3-r-7-direct}
except that the first $d$ sets of 12 triples of the base block set $B$ are $0^\mathrm{opp}$, $1^\mathrm{opp}$, \dots, $(d-1)^\mathrm{opp}$.

For each of the remaining values of $r$, create a random connected $(9,5^+)$-graph with $2r + 10$ vertices.
On the neighbours of each point $x$, construct a Steiner triple system of order 9 to form the blocks of $x^\mathrm{opp}$.
Create the remaining blocks of the geometry by hill climbing.
The process is readily automated.


{\boldmath $\adfPENT(3,57,9)$}, $d = 4$:
{21, 29, 35;  21, 54, 103;  21, 55, 115;  21, 70, 113;  29, 54, 115;  29, 55, 113;  29, 70, 103;  35, 54, 113;  35, 55, 103;  35, 70, 115;  54, 55, 70;  103, 113, 115;  10, 12, 19;  10, 22, 96;  10, 70, 107;  10, 90, 104;  12, 22, 107;  12, 70, 104;  12, 90, 96;  19, 22, 104;  19, 70, 96;  19, 90, 107;  22, 70, 90;  96, 104, 107;  23, 31, 37;  23, 56, 105;  23, 57, 117;  23, 72, 115;  31, 56, 117;  31, 57, 115;  31, 72, 105;  37, 56, 115;  37, 57, 105;  37, 72, 117;  56, 57, 72;  105, 115, 117;  12, 14, 21;  12, 24, 98;  12, 72, 109;  12, 92, 106;  14, 24, 109;  14, 72, 106;  14, 92, 98;  21, 24, 106;  21, 72, 98;  21, 92, 109;  24, 72, 92;  98, 106, 109;  0, 4, 67;  0, 5, 22;  0, 13, 36;  0, 15, 117;  0, 18, 31;  0, 19, 81;  0, 23, 27;  0, 24, 52;  0, 25, 77;  0, 41, 69;  0, 51, 123;  0, 57, 102;  0, 62, 87;  0, 65, 119;  0, 71, 93;  0, 79, 107;  0, 106, 111;  1, 2, 17;  1, 5, 98;  1, 6, 71;  1, 14, 101;  1, 54, 106;  1, 58, 62;  1, 74, 102;  2, 26, 95;  2, 38, 79;  2, 59, 83;  2, 103, 119}

{\boldmath $\adfPENT(3,58,9)$}, $d = 6$:
{23, 34, 49;  23, 51, 92;  23, 59, 111;  23, 63, 105;  34, 51, 111;  34, 59, 105;  34, 63, 92;  49, 51, 105;  49, 59, 92;  49, 63, 111;  51, 59, 63;  92, 105, 111;  16, 17, 22;  16, 64, 78;  16, 68, 111;  16, 76, 104;  17, 64, 111;  17, 68, 104;  17, 76, 78;  22, 64, 104;  22, 68, 78;  22, 76, 111;  64, 68, 76;  78, 104, 111;  25, 36, 51;  25, 53, 94;  25, 61, 113;  25, 65, 107;  36, 53, 113;  36, 61, 107;  36, 65, 94;  51, 53, 107;  51, 61, 94;  51, 65, 113;  53, 61, 65;  94, 107, 113;  18, 19, 24;  18, 66, 80;  18, 70, 113;  18, 78, 106;  19, 66, 113;  19, 70, 106;  19, 78, 80;  24, 66, 106;  24, 70, 80;  24, 78, 113;  66, 70, 78;  80, 106, 113;  27, 38, 53;  27, 55, 96;  27, 63, 115;  27, 67, 109;  38, 55, 115;  38, 63, 109;  38, 67, 96;  53, 55, 109;  53, 63, 96;  53, 67, 115;  55, 63, 67;  96, 109, 115;  20, 21, 26;  20, 68, 82;  20, 72, 115;  20, 80, 108;  21, 68, 115;  21, 72, 108;  21, 80, 82;  26, 68, 108;  26, 72, 82;  26, 80, 115;  68, 72, 80;  82, 108, 115;  0, 3, 20;  0, 5, 101;  0, 9, 113;  0, 11, 117;  0, 16, 27;  0, 18, 55;  0, 21, 102;  0, 22, 91;  0, 30, 61;  0, 32, 69;  0, 41, 50;  0, 53, 73;  0, 76, 107;  0, 81, 94;  0, 87, 106;  0, 99, 123;  0, 103, 125;  0, 104, 109;  0, 110, 119;  1, 2, 26;  1, 4, 100;  1, 8, 77;  1, 10, 97;  1, 14, 25;  1, 19, 46;  1, 20, 93;  1, 21, 106;  1, 35, 69;  1, 51, 124;  1, 59, 82;  1, 74, 105;  1, 86, 118;  2, 5, 32;  2, 20, 47;  2, 22, 83;  2, 23, 125;  2, 52, 57;  2, 63, 93;  2, 89, 106;  2, 105, 112;  3, 4, 28;  3, 21, 71;  4, 22, 95;  4, 41, 59}
\end{proof}
%
\begin{theorem}
\label{thm:PENT-3-r-9-constructed}
There exists a $\adfPENT(3,r,9)$ with deficiency graph of girth at least $5$ if $r \equiv 0 \textrm{~or~} 1 \adfmod{3}$ and $r \ge 57$.
\end{theorem}
\begin{proof}
The main argument is similar to that of Theorem~\ref{thm:PENT-3-r-7-constructed}.

Let $S_9 = \{s: 57 \le s \le 118,\; s \equiv 0 \textrm{~or~} 1 \adfmod{3}\}$ and observe that
this set covers the admissible residues modulo 63.
By Lemma~\ref{lem:PENT-3-r-9-direct}, there exists a $\adfPENT(3,s,9)$ for $s \in S_9$, and
it suffices to address only $r \ge 246$.

For each $s \in S_9$, take $u$ copies of a $\adfPENT(3,58,9)$, which has 126 points, and one copy of a $\adfPENT(3,s,9)$.
By Lemma~\ref{lem:3-GDD-existence}, since $2 s + 10 \le 126\cdot 2$,
there exists a 3-GDD of type $126^u (2s + 10)^1$ for $u = 3$, 4, \dots.
Now use Theorem~\ref{thm:GDD-basic} to construct a $\adfPENT(3,63u + s,9)$ for $u \ge 3$.

Hence there exists a $\adfPENT(3,r,9)$ for admissible $r \ge 3 \cdot 63 + 57 = 246$.
\end{proof}
It is clear that one can continue in a similar manner until one's computers run out of steam.
With the availability of sufficient resources,
`starter systems' can be obtained to prove similar theorems for $w > 9$ that are valid as orders of Steiner triple systems.
We conclude this section by dealing with $w \in$ \{13, 15, 19, 21, 25, 27, 31, 33\}.
%
\begin{lemma}
\label{lem:PENT-3-r-w-direct}
There exists a $\adfPENT(3,r,13)$ with connected deficiency graph of girth at least $5$ for $168 \le r \le 698$,
$r \equiv 0 \textrm{~or~} 2 \adfmod{3}$.

There exists a $\adfPENT(3,r,15)$ with connected deficiency graph of girth at least $5$ for $249 \le r \le 1021$,
$r \equiv 0 \textrm{~or~} 1 \adfmod{3}$.

There exists a $\adfPENT(3,r,19)$ with connected deficiency graph of girth at least $5$ for $479 \le r \le 1455$,
$r \equiv 3 \textrm{~or~} 5 \adfmod{6}$.

There exists a $\adfPENT(3,r,21)$ with connected deficiency graph of girth at least $5$ for $654 \le r \le 1990$,
$r \equiv 0 \textrm{~or~} 4 \adfmod{6}$.

There exists a $\adfPENT(3,r,w)$ with connected deficiency graph of girth at least $5$ for
\begin{align*}
(r,w) \in \{&(1212, 25), (1214, 25), (1220, 25), (1226, 25), \\
            &(1389, 27), (1393, 27), (1399, 27), (1405, 27), \\
            &(2133, 31), (2135, 31), (2141, 31), (2147, 31), \\
            &(2550, 33), (2554, 33), (2560, 33), (2566, 33)\}.
\end{align*}
\end{lemma}
\begin{proof}
Let $w \in$ \{13, 15, 19, 21, 25, 27, 31, 33\}.
For each $r$, create a random connected $(w,5^+)$-graph with $v = 2r + w + 1$ vertices.
On the neighbours of each point $x$
construct a Steiner triple system of order $w$ to form the blocks of $x^\mathrm{opp}$.
Create the remaining blocks of the geometry by hill climbing.
We give some rather bulky examples (including all sixteen for $w \in$ \{25, 27, 31, 33\}) in the Appendix.
Alternatively, you may email the first
author for the entire collection, including those of Lemmas~\ref{lem:PENT-3-r-7-direct} and \ref{lem:PENT-3-r-9-direct}.
We used the cyclic STS(13) and the projective STS(15), \#1 in \cite[Table 5.8]{ColbourRosa1999}.

Observe that when $w \ge 19$ we restrict $r$ to values where
the number of points in the corresponding $\adfPENT(3, r, w)$ is congruent to 2 modulo 4.
With this constraint we are able to exploit an automorphism of the form
$x \mapsto x + 2 \adfmod{v}$ or $x \mapsto x + 6 \adfmod{v}$.
\end{proof}
%
\begin{theorem}
\label{thm:PENT-3-r-13-15-constructed}
There exists a $\adfPENT(3,r,13)$ with deficiency graph of girth at least $5$ if $r \equiv 0 \textrm{~or~} 2 \adfmod{3}$ and $r \ge 168$.

There exists a $\adfPENT(3,r,15)$ with deficiency graph of girth at least $5$ if $r \equiv 0 \textrm{~or~} 1 \adfmod{3}$ and $r \ge 249$.
\end{theorem}
\begin{proof}
Let $w \in \{13, 15\}$,
\begin{align*}
S_{13} &= \{s: 168 \le s \le 344,\; s \equiv 0 \textrm{~or~} 2 \adfmod{3}\}, \\
S_{15} &= \{s: 249 \le s \le 505,\; s \equiv 0 \textrm{~or~} 1 \adfmod{3}\},
\end{align*}
$a_{13} = 170$, $a_{15} = 250$, and $m_w = a_w + (w + 1)/2$.

By Lemma~\ref{lem:PENT-3-r-w-direct}, there exists a $\adfPENT(3,s,w)$ for all $s \in S_{w}$.
Also the elements of $S_w$ cover the admissible residue classes modulo $m_w$,
and $S_w$ consists of all admissible values of $s$ that occur in the specified range.

Let $s \in S_w$.
Take $u$ copies of a $\adfPENT(3,a_w,w)$, which has $2a_w + w + 1$ points, and one copy of a $\adfPENT(3,s,w)$.
By Lemma~\ref{lem:3-GDD-existence}, since $2a_w + w + 1$ is a multiple of 3 and $2 s + w + 1 \le 2(2a_w + w + 1)$,
there exists a 3-GDD of type $(2a_w + w + 1)^u (2s + w + 1)^1$ for $u = 3$, 4, \dots.
Now use Theorem~\ref{thm:GDD-basic} to construct a $\adfPENT(3, m_w u + s, w)$ for $u \ge 3$.

Hence there exists a $\adfPENT(3, r, w)$ for all admissible $r \ge 3 \cdot m_w  + \min(S_w)$.
The remaining $r$ to complete the proof are given by Lemma~\ref{lem:PENT-3-r-w-direct}.
\end{proof}
%
\begin{theorem}
\label{thm:PENT-3-r-19-21-constructed}
There exists a $\adfPENT(3,r,19)$ with deficiency graph of girth at least $5$ if $r \equiv 0 \textrm{~or~} 2 \adfmod{3}$ and $r \ge 3411$.

There exists a $\adfPENT(3,r,21)$ with deficiency graph of girth at least $5$ if $r \equiv 0 \textrm{~or~} 1 \adfmod{3}$ and $r \ge 4666$.
\end{theorem}
\begin{proof}
Let $w \in \{19, 21\}$,
\begin{align*}
S_{19} &= \{s: 479 \le s \le 1455,\; s \equiv 3 \textrm{~or~} 5 \adfmod{6}\}, \\
S_{21} &= \{s: 654 \le s \le 1990,\; s \equiv 0 \textrm{~or~} 4 \adfmod{6}\},
\end{align*}
$a_{19} = 479$, $a_{21} = 658$ and $m_w = a_w + (w + 1)/2$.

By Lemma~\ref{lem:PENT-3-r-w-direct}, there exists a $\adfPENT(3,s,w)$ for all $s \in S_{w}$.
Also the elements of $S_w$ cover the admissible residue classes modulo $m_w$, but
note that $S_w$ contains approximately half of the admissible values of $s$ in the specified range.

Let $s \in S_w$.
Take $u$ copies of a $\adfPENT(3,a_w,w)$, which has $2a_w + w + 1$ points, and one copy of a $\adfPENT(3,s,w)$.
By Lemma~\ref{lem:3-GDD-existence}, since $2a_w + w + 1$ is a multiple of 3 and $2 s + w + 1 \le 3(2a_w + w + 1)$,
there exists a 3-GDD of type $(2a_w + w + 1)^u (2s + w + 1)^1$ for $u = 4$, 5, \dots.
Now use Theorem~\ref{thm:GDD-basic} to construct a $\adfPENT(3, m_w u + s, w)$ for $u \ge 4$.

Hence there exists a $\adfPENT(3, r, w)$ for all admissible $r \ge 4 \cdot m_w  + \max(S_w)$.
\end{proof}
%
\begin{theorem}
\label{thm:PENT-3-r-w-constructed-large}
There exist generalized pentagonal geometries $\adfPENT(3, r, w)$ with deficiency graph of girth at least $5$
for $w \in \{25, 27, 31, 33\}$ and all sufficiently large admissible $r$.
\end{theorem}
\begin{proof}
Let $w \in \{25, 27, 31, 33\}$. Let
$$\begin{array}{ll}
  r_0 = 1212,~~r_1 = 1214,~~r_2 = 1220,~~ r_3 = 1226 & \text{~if~} w = 25, \\
  r_0 = 1389,~~r_1 = 1393,~~r_2 = 1399,~~ r_3 = 1405 & \text{~if~} w = 27, \\
  r_0 = 2133,~~r_1 = 2135,~~r_2 = 3241,~~ r_3 = 2147 & \text{~if~} w = 31, \\
  r_0 = 2550,~~r_1 = 2554,~~r_2 = 2560,~~ r_3 = 2566 & \text{~if~} w = 33.
\end{array}$$
Let $v_i = 2 r_i + w + 1$, $i = 1$, 2, 3, 4.
Then by Lemma~\ref{lem:PENT-3-r-w-direct} there exist a $\adfPENT(3, r, w)$ for
$r = r_0$, $r_1$, $r_2$, $r_3$.
Also $r_0 \equiv v_1 \equiv v_2 \equiv v_3 \equiv 0 \adfmod{3}$.
Let $r_* \in \{r_0, r_3\}$ and $v_* = 2r_* + w + 1$.

Apply Theorem~\ref{thm:GDD-basic} twice, first with
a $\adfPENT(3, r_*, w)$,
$t$ copies of a $\adfPENT(3, r_1, w)$
and a 3-GDD of type $v_1^t v_*^1$ to obtain a
$$\adfPENT(3, v_1 t/2 + r_*, w)~~ \text{~for~} t \ge 3,$$
and then with
$u$ copies of a $\adfPENT(3, r_2, w)$
and a 3-GDD of type $v_2^u (v_1 t + v_*)^1$ to obtain a
\begin{equation}\label{eqn:PENT-3-r-w-constructed}
\adfPENT(3, v_2 u/2 + v_1 t/2 + r_*, w)~~ \text{~for~} u \ge t + 2,~~ t \ge 3.
\end{equation}
To confirm the existence of the two 3-GDDs by Lemma~\ref{lem:3-GDD-existence}, we have
$v_1 \equiv  v_2 \equiv 0 \adfmod{6}$, $v_* \equiv 0 \adfmod{2}$,
$v_* \le 2 v_1 \le v_1(t - 1)$ and
$v_1 t + v_* \le (v_2 - 12)t + v_2 + 32 \le v_2(t + 1) \le v_2(u - 1)$.

Observe that in (\ref{eqn:PENT-3-r-w-constructed}) the greatest common divisor of the coefficients of $u$ and $t$ is 3, and
that the two options for $r_*$ belong to the two admissible residue classes modulo 3.
Hence, by an argument involving the Chinese Remainder Theorem similar to that in the proof of Theorem~\ref{thm:PENT-4-r-13-constructed}, below,
the range of $v_2 u/2 + v_1 t/2 + r_*$ in (\ref{eqn:PENT-3-r-w-constructed})
covers all sufficiently large admissible values.
\end{proof}


\section{Block size 4}\label{sec:Block size 4}
By Lemma~\ref{lem:v-b}, a pentagonal geometry $\adfPENT(4,r,w)$ has $3r + w + 1$ points and $r(3r + w + 1)/4$ lines.
Furthermore, a Steiner system $S(2,4,w)$ with positive $w$ exists if and only if $w \equiv 1 \textrm{~or~} 4 \adfmod{12}$.
Therefore $r$ is admissible if and only if $r \equiv 0 \textrm{~or~} w + 1 \adfmod{4}$.
In contrast to the case $k = 3$,
examples of $\adfPENT(4,r,w)$ with $w > 4$ and pairwise disjoint opposite designs appear to be quite difficult to find.
The following lemma represents what we have obtained.
%
\begin{lemma}
\label{lem:PENT-4-r-13-direct}
There exists a $\adfPENT(4,r,13)$ with connected deficiency graph of girth $5$ if
$r \in$ \{$112$, $116$, $120$, $124$, $128$, $132$, $136$, $140$, $144$, $148$, $152$, $156$, $160$, $164$\}.
\end{lemma}
\begin{proof}
In each case we specify $d = 2$ followed by $2r$ numbers.
To obtain the geometry, gather the numbers into quadruples to form a set of $r/2$ base blocks, $B$,
which are developed into the line set of the geometry by the mapping $x \mapsto x + 2 \adfmod{3r + 14}$.
The first two sets of 13 quadruples of $B$ are $0^\mathrm{opp}$ and $1^\mathrm{opp}$.
The opposite designs are projective planes of order 3.
The deficiency graphs have girth 5.


{\boldmath $\adfPENT(4,112,13)$}, $d = 2$:
{21, 41, 49, 297;  21, 140, 210, 319;  21, 149, 293, 335;  21, 163, 287, 337;  41, 140, 293, 337;  41, 149, 287, 319;  41, 163, 210, 335;  49, 140, 287, 335;  49, 149, 210, 337;  49, 163, 293, 319;  140, 149, 163, 297;  210, 287, 293, 297;  297, 319, 335, 337;  14, 16, 32, 293;  14, 54, 64, 302;  14, 58, 202, 310;  14, 59, 188, 330;  16, 54, 202, 330;  16, 58, 188, 302;  16, 59, 64, 310;  32, 54, 188, 310;  32, 58, 64, 330;  32, 59, 202, 302;  54, 58, 59, 293;  64, 188, 202, 293;  293, 302, 310, 330;  110, 249, 50, 261;  206, 137, 183, 48;  247, 154, 277, 34;  132, 347, 323, 292;  34, 58, 37, 189;  167, 71, 106, 340;  196, 346, 269, 11;  0, 7, 46, 291;  0, 13, 133, 136;  0, 11, 137, 264;  0, 17, 198, 217;  0, 25, 110, 185;  0, 68, 209, 269;  0, 51, 119, 126;  0, 58, 187, 305;  0, 63, 82, 323;  0, 35, 76, 145;  0, 33, 88, 307;  0, 47, 205, 254;  0, 115, 166, 325;  0, 29, 117, 284;  0, 92, 203, 349;  0, 39, 118, 175;  0, 113, 146, 267;  0, 90, 193, 263;  0, 53, 154, 321;  0, 85, 171, 182;  0, 65, 132, 233;  0, 12, 71, 237;  0, 30, 67, 285}

{\boldmath $\adfPENT(4,116,13)$}, $d = 2$:
{9, 20, 39, 277;  9, 55, 141, 311;  9, 63, 267, 329;  9, 77, 147, 342;  20, 55, 267, 342;  20, 63, 147, 311;  20, 77, 141, 329;  39, 55, 147, 329;  39, 63, 141, 342;  39, 77, 267, 311;  55, 63, 77, 277;  141, 147, 267, 277;  277, 311, 329, 342;  34, 52, 57, 307;  34, 86, 222, 308;  34, 96, 300, 324;  34, 216, 286, 354;  52, 86, 300, 354;  52, 96, 286, 308;  52, 216, 222, 324;  57, 86, 286, 324;  57, 96, 222, 354;  57, 216, 300, 308;  86, 96, 216, 307;  222, 286, 300, 307;  307, 308, 324, 354;  217, 15, 206, 64;  84, 65, 236, 39;  159, 303, 66, 34;  206, 105, 304, 221;  148, 79, 197, 264;  268, 309, 305, 82;  3, 294, 83, 220;  81, 8, 237, 192;  6, 107, 268, 9;  206, 341, 38, 254;  159, 314, 326, 101;  350, 231, 44, 307;  281, 159, 127, 168;  0, 1, 28, 281;  0, 61, 66, 239;  0, 58, 175, 341;  0, 80, 189, 289;  0, 47, 144, 359;  0, 112, 235, 271;  0, 31, 199, 312;  0, 4, 217, 245;  0, 33, 53, 326;  0, 67, 219, 244;  0, 99, 160, 275;  0, 17, 233, 240;  0, 51, 154, 347;  0, 25, 27, 82;  0, 129, 169, 206;  0, 13, 76, 155;  0, 131, 143, 166;  0, 26, 183, 231;  0, 2, 107, 181}

{\boldmath $\adfPENT(4,120,13)$}, $d = 2$:
{11, 18, 31, 309;  11, 97, 263, 325;  11, 109, 295, 356;  11, 261, 269, 363;  18, 97, 295, 363;  18, 109, 269, 325;  18, 261, 263, 356;  31, 97, 269, 356;  31, 109, 263, 363;  31, 261, 295, 325;  97, 109, 261, 309;  263, 269, 295, 309;  309, 325, 356, 363;  12, 50, 66, 278;  12, 71, 112, 305;  12, 80, 266, 344;  12, 106, 114, 364;  50, 71, 266, 364;  50, 80, 114, 305;  50, 106, 112, 344;  66, 71, 114, 344;  66, 80, 112, 364;  66, 106, 266, 305;  71, 80, 106, 278;  112, 114, 266, 278;  278, 305, 344, 364;  236, 99, 308, 47;  55, 46, 185, 288;  162, 41, 371, 246;  108, 39, 339, 314;  290, 145, 323, 262;  60, 295, 152, 101;  247, 192, 189, 58;  310, 300, 367, 185;  160, 83, 343, 340;  132, 355, 4, 115;  32, 210, 155, 348;  9, 116, 235, 220;  0, 1, 19, 322;  0, 4, 107, 304;  0, 17, 145, 170;  0, 23, 161, 262;  0, 35, 45, 324;  0, 65, 137, 321;  0, 24, 171, 341;  0, 44, 133, 217;  0, 51, 192, 275;  0, 121, 140, 289;  0, 43, 175, 190;  0, 75, 201, 292;  0, 63, 219, 224;  0, 37, 87, 256;  0, 47, 151, 226;  0, 73, 77, 244;  0, 126, 285, 373;  0, 114, 301, 329;  0, 105, 197, 268;  0, 88, 241, 353;  0, 69, 337, 361;  0, 93, 129, 156}

{\boldmath $\adfPENT(4,124,13)$}, $d = 2$:
{19, 37, 82, 285;  19, 151, 211, 304;  19, 205, 247, 311;  19, 207, 239, 355;  37, 151, 247, 355;  37, 205, 239, 304;  37, 207, 211, 311;  82, 151, 239, 311;  82, 205, 211, 355;  82, 207, 247, 304;  151, 205, 207, 285;  211, 239, 247, 285;  285, 304, 311, 355;  32, 76, 99, 236;  32, 102, 176, 289;  32, 140, 182, 350;  32, 148, 180, 368;  76, 102, 182, 368;  76, 140, 180, 289;  76, 148, 176, 350;  99, 102, 180, 350;  99, 140, 176, 368;  99, 148, 182, 289;  102, 140, 148, 236;  176, 180, 182, 236;  236, 289, 350, 368;  230, 363, 94, 79;  258, 335, 95, 306;  252, 63, 80, 93;  30, 249, 25, 54;  30, 281, 223, 14;  92, 336, 5, 16;  252, 209, 3, 62;  79, 38, 13, 3;  380, 209, 78, 69;  191, 312, 65, 214;  374, 168, 279, 23;  15, 237, 308, 12;  145, 332, 167, 214;  260, 233, 2, 1;  61, 365, 326, 164;  293, 214, 204, 231;  0, 1, 143, 373;  0, 11, 279, 300;  0, 12, 179, 331;  0, 47, 95, 264;  0, 45, 52, 178;  0, 65, 156, 243;  0, 22, 281, 297;  0, 5, 57, 130;  0, 15, 277, 328;  0, 49, 246, 349;  0, 14, 173, 185;  0, 30, 63, 191;  0, 43, 153, 284;  0, 59, 181, 353;  0, 9, 99, 240;  0, 25, 110, 234;  0, 21, 271, 324;  0, 71, 91, 294;  0, 31, 115, 154;  0, 20, 81, 105}

{\boldmath $\adfPENT(4,128,13)$}, $d = 2$:
{341, 43, 86, 117;  341, 71, 312, 349;  341, 351, 381, 165;  341, 97, 121, 55;  43, 71, 381, 55;  43, 351, 121, 349;  43, 97, 312, 165;  86, 71, 121, 165;  86, 351, 312, 55;  86, 97, 381, 349;  71, 351, 97, 117;  312, 121, 381, 117;  117, 349, 165, 55;  302, 18, 58, 234;  302, 50, 282, 328;  302, 344, 356, 278;  302, 151, 249, 48;  18, 50, 356, 48;  18, 344, 249, 328;  18, 151, 282, 278;  58, 50, 249, 278;  58, 344, 282, 48;  58, 151, 356, 328;  50, 344, 151, 234;  282, 249, 356, 234;  234, 328, 278, 48;  55, 345, 156, 118;  291, 202, 286, 14;  251, 178, 55, 195;  184, 97, 364, 308;  291, 131, 136, 130;  79, 220, 143, 144;  287, 176, 274, 15;  207, 166, 230, 131;  299, 148, 310, 357;  94, 57, 265, 116;  192, 19, 269, 174;  23, 180, 286, 1;  140, 3, 239, 10;  293, 151, 178, 260;  112, 203, 51, 254;  0, 3, 248, 395;  0, 9, 197, 238;  0, 14, 204, 249;  0, 7, 233, 318;  0, 21, 59, 140;  0, 19, 132, 213;  0, 108, 245, 331;  0, 65, 158, 283;  0, 173, 192, 373;  0, 23, 242, 385;  0, 107, 159, 239;  0, 83, 200, 323;  0, 53, 177, 259;  0, 105, 175, 196;  0, 63, 250, 325;  0, 61, 329, 343;  0, 58, 167, 185;  0, 52, 269, 353;  0, 27, 116, 152;  0, 57, 96, 163;  0, 25, 136, 389;  0, 34, 85, 327;  0, 49, 145, 328}

{\boldmath $\adfPENT(4,132,13)$}, $d = 2$:
{5, 23, 75, 289;  5, 89, 187, 377;  5, 129, 278, 399;  5, 132, 247, 403;  23, 89, 278, 403;  23, 129, 247, 377;  23, 132, 187, 399;  75, 89, 247, 399;  75, 129, 187, 403;  75, 132, 278, 377;  89, 129, 132, 289;  187, 247, 278, 289;  289, 377, 399, 403;  8, 12, 34, 322;  8, 105, 224, 336;  8, 122, 307, 388;  8, 164, 282, 406;  12, 105, 307, 406;  12, 122, 282, 336;  12, 164, 224, 388;  34, 105, 282, 388;  34, 122, 224, 406;  34, 164, 307, 336;  105, 122, 164, 322;  224, 282, 307, 322;  322, 336, 388, 406;  77, 262, 114, 218;  132, 298, 345, 111;  80, 401, 340, 183;  154, 64, 0, 81;  340, 166, 119, 375;  69, 392, 347, 120;  400, 233, 29, 380;  38, 249, 313, 323;  55, 229, 264, 188;  63, 76, 68, 113;  27, 291, 54, 132;  74, 392, 221, 159;  397, 182, 40, 259;  170, 15, 17, 277;  84, 225, 290, 156;  93, 86, 262, 312;  0, 1, 29, 408;  0, 10, 63, 205;  0, 6, 15, 181;  0, 19, 27, 94;  0, 13, 32, 197;  0, 21, 113, 374;  0, 28, 95, 173;  0, 24, 341, 347;  0, 133, 153, 230;  0, 111, 135, 218;  0, 46, 119, 151;  0, 116, 279, 315;  0, 101, 169, 217;  0, 43, 123, 362;  0, 33, 161, 240;  0, 62, 297, 369;  0, 231, 325, 401;  0, 249, 293, 339;  0, 120, 251, 385;  0, 49, 281, 342;  0, 65, 188, 355;  0, 139, 178, 361;  0, 79, 80, 309;  0, 51, 74, 202}

{\boldmath $\adfPENT(4,136,13)$}, $d = 2$:
{343, 209, 397, 416;  343, 325, 1, 301;  343, 287, 119, 6;  343, 389, 49, 51;  209, 325, 119, 51;  209, 287, 49, 301;  209, 389, 1, 6;  397, 325, 49, 6;  397, 287, 1, 51;  397, 389, 119, 301;  325, 287, 389, 416;  1, 49, 119, 416;  416, 301, 6, 51;  304, 98, 374, 413;  304, 214, 0, 136;  304, 122, 80, 11;  304, 372, 26, 34;  98, 214, 80, 34;  98, 122, 26, 136;  98, 372, 0, 11;  374, 214, 26, 11;  374, 122, 0, 34;  374, 372, 80, 136;  214, 122, 372, 413;  0, 26, 80, 413;  413, 136, 11, 34;  334, 59, 268, 339;  371, 50, 228, 145;  56, 76, 339, 213;  148, 89, 288, 321;  180, 220, 297, 30;  295, 168, 191, 131;  207, 102, 132, 256;  28, 312, 337, 65;  366, 29, 240, 51;  24, 130, 295, 229;  52, 233, 23, 166;  87, 145, 240, 292;  365, 138, 334, 57;  80, 390, 79, 335;  400, 215, 129, 206;  391, 408, 101, 162;  186, 369, 306, 123;  234, 5, 296, 202;  128, 5, 340, 57;  293, 44, 289, 333;  0, 3, 290, 305;  0, 10, 27, 89;  0, 29, 200, 204;  0, 19, 60, 318;  0, 13, 187, 280;  0, 22, 265, 377;  0, 86, 317, 347;  0, 103, 156, 323;  0, 73, 260, 401;  0, 16, 149, 185;  0, 58, 179, 371;  0, 21, 221, 415;  0, 163, 195, 303;  0, 69, 81, 314;  0, 161, 174, 365;  0, 257, 273, 411;  0, 36, 83, 207;  0, 28, 93, 230;  0, 35, 256, 357;  0, 91, 97, 253;  0, 53, 159, 322;  0, 44, 241, 419}

{\boldmath $\adfPENT(4,140,13)$}, $d = 2$:
{8, 9, 109, 379;  8, 159, 295, 383;  8, 203, 327, 401;  8, 213, 321, 426;  9, 159, 327, 426;  9, 203, 321, 383;  9, 213, 295, 401;  109, 159, 321, 401;  109, 203, 295, 426;  109, 213, 327, 383;  159, 203, 213, 379;  295, 321, 327, 379;  379, 383, 401, 426;  34, 52, 56, 276;  34, 79, 140, 326;  34, 108, 232, 357;  34, 114, 222, 426;  52, 79, 232, 426;  52, 108, 222, 326;  52, 114, 140, 357;  56, 79, 222, 357;  56, 108, 140, 426;  56, 114, 232, 326;  79, 108, 114, 276;  140, 222, 232, 276;  276, 326, 357, 426;  388, 155, 59, 52;  364, 352, 351, 20;  268, 192, 185, 420;  350, 385, 70, 118;  140, 142, 75, 375;  281, 382, 301, 421;  203, 398, 57, 142;  389, 403, 312, 8;  153, 304, 170, 1;  239, 404, 366, 43;  14, 240, 369, 61;  260, 3, 79, 403;  421, 424, 58, 242;  62, 201, 48, 333;  272, 257, 160, 85;  293, 338, 314, 17;  244, 388, 147, 0;  101, 34, 191, 64;  126, 267, 405, 221;  172, 94, 345, 357;  150, 339, 369, 350;  0, 3, 5, 96;  0, 13, 40, 115;  0, 11, 172, 347;  0, 15, 34, 197;  0, 57, 276, 361;  0, 25, 63, 120;  0, 59, 171, 348;  0, 43, 66, 275;  0, 41, 72, 241;  0, 128, 261, 309;  0, 21, 149, 165;  0, 28, 183, 324;  0, 71, 79, 107;  0, 83, 289, 294;  0, 29, 53, 398;  0, 20, 119, 345;  0, 69, 146, 278;  0, 49, 121, 238;  0, 126, 385, 425;  0, 70, 231, 297;  0, 51, 174, 229;  0, 123, 191, 277;  0, 33, 207, 397}

{\boldmath $\adfPENT(4,144,13)$}, $d = 2$:
{73, 93, 99, 267;  73, 101, 145, 336;  73, 110, 225, 431;  73, 131, 155, 443;  93, 101, 225, 443;  93, 110, 155, 336;  93, 131, 145, 431;  99, 101, 155, 431;  99, 110, 145, 443;  99, 131, 225, 336;  101, 110, 131, 267;  145, 155, 225, 267;  267, 336, 431, 443;  4, 16, 129, 346;  4, 180, 302, 348;  4, 222, 319, 354;  4, 292, 316, 374;  16, 180, 319, 374;  16, 222, 316, 348;  16, 292, 302, 354;  129, 180, 316, 354;  129, 222, 302, 374;  129, 292, 319, 348;  180, 222, 292, 346;  302, 316, 319, 346;  346, 348, 354, 374;  382, 283, 296, 86;  161, 7, 251, 168;  178, 307, 160, 239;  222, 218, 158, 29;  135, 72, 87, 364;  213, 416, 333, 158;  220, 297, 48, 67;  108, 18, 141, 67;  129, 311, 116, 312;  11, 54, 372, 29;  160, 238, 99, 373;  339, 402, 335, 275;  219, 428, 250, 272;  365, 3, 42, 308;  354, 281, 120, 189;  321, 59, 148, 378;  134, 27, 254, 293;  103, 357, 206, 317;  7, 66, 206, 257;  34, 388, 1, 338;  123, 321, 204, 37;  84, 48, 409, 132;  357, 66, 320, 171;  198, 164, 412, 317;  0, 1, 242, 372;  0, 5, 106, 341;  0, 16, 41, 179;  0, 98, 269, 399;  0, 7, 246, 313;  0, 43, 66, 289;  0, 23, 89, 224;  0, 53, 87, 238;  0, 91, 182, 359;  0, 68, 299, 397;  0, 71, 121, 186;  0, 167, 189, 425;  0, 9, 118, 305;  0, 199, 215, 441;  0, 29, 251, 369;  0, 11, 47, 308;  0, 40, 184, 327;  0, 100, 233, 375;  0, 31, 239, 284;  0, 65, 202, 367;  0, 81, 181, 256;  0, 137, 349, 427}

{\boldmath $\adfPENT(4,148,13)$}, $d = 2$:
{79, 41, 147, 259;  79, 96, 49, 283;  79, 111, 362, 277;  79, 417, 415, 427;  41, 96, 362, 427;  41, 111, 415, 283;  41, 417, 49, 277;  147, 96, 415, 277;  147, 111, 49, 427;  147, 417, 362, 283;  96, 111, 417, 259;  49, 415, 362, 259;  259, 283, 277, 427;  44, 32, 182, 200;  44, 101, 42, 348;  44, 176, 359, 312;  44, 410, 380, 418;  32, 101, 359, 418;  32, 176, 380, 348;  32, 410, 42, 312;  182, 101, 380, 312;  182, 176, 42, 418;  182, 410, 359, 348;  101, 176, 410, 200;  42, 380, 359, 200;  200, 348, 312, 418;  3, 322, 198, 187;  396, 52, 129, 223;  111, 302, 229, 396;  244, 284, 265, 156;  256, 242, 15, 93;  330, 5, 54, 415;  194, 217, 340, 345;  385, 352, 410, 394;  435, 183, 448, 391;  23, 180, 374, 111;  431, 33, 317, 226;  44, 18, 426, 317;  387, 146, 241, 227;  404, 296, 413, 178;  418, 199, 113, 314;  126, 145, 91, 78;  119, 284, 123, 306;  446, 276, 237, 301;  433, 2, 169, 68;  102, 215, 189, 139;  285, 440, 338, 227;  196, 119, 385, 224;  176, 130, 251, 428;  0, 1, 414, 443;  0, 3, 43, 372;  0, 4, 212, 219;  0, 17, 20, 196;  0, 27, 342, 457;  0, 31, 201, 406;  0, 54, 209, 421;  0, 71, 78, 333;  0, 99, 296, 453;  0, 197, 305, 429;  0, 103, 174, 397;  0, 39, 61, 294;  0, 123, 169, 256;  0, 93, 126, 200;  0, 63, 165, 394;  0, 35, 237, 358;  0, 97, 341, 393;  0, 131, 214, 413;  0, 184, 369, 435;  0, 105, 125, 221;  0, 56, 289, 363;  0, 60, 287, 329;  0, 213, 339, 395;  0, 73, 186, 220;  0, 119, 357, 391}

{\boldmath $\adfPENT(4,152,13)$}, $d = 2$:
{26, 105, 125, 365;  26, 127, 259, 389;  26, 165, 307, 419;  26, 227, 293, 444;  105, 127, 307, 444;  105, 165, 293, 389;  105, 227, 259, 419;  125, 127, 293, 419;  125, 165, 259, 444;  125, 227, 307, 389;  127, 165, 227, 365;  259, 293, 307, 365;  365, 389, 419, 444;  52, 82, 106, 344;  52, 111, 212, 346;  52, 164, 306, 361;  52, 178, 244, 366;  82, 111, 306, 366;  82, 164, 244, 346;  82, 178, 212, 361;  106, 111, 244, 361;  106, 164, 212, 366;  106, 178, 306, 346;  111, 164, 178, 344;  212, 244, 306, 344;  344, 346, 361, 366;  38, 20, 251, 225;  383, 270, 341, 406;  426, 208, 227, 318;  135, 153, 198, 128;  134, 447, 411, 362;  404, 8, 405, 417;  398, 352, 208, 315;  450, 361, 120, 143;  256, 417, 272, 152;  12, 469, 15, 365;  251, 58, 284, 280;  249, 199, 20, 76;  187, 230, 339, 258;  70, 394, 273, 1;  290, 149, 155, 380;  23, 159, 324, 150;  255, 456, 180, 27;  153, 100, 24, 349;  70, 184, 127, 331;  372, 257, 205, 72;  231, 204, 319, 20;  80, 88, 44, 300;  444, 347, 455, 61;  122, 443, 128, 19;  223, 59, 178, 188;  133, 72, 431, 394;  0, 21, 31, 428;  0, 33, 118, 157;  0, 43, 47, 420;  0, 51, 346, 429;  0, 65, 92, 181;  0, 12, 171, 318;  0, 98, 371, 455;  0, 95, 274, 395;  0, 64, 167, 311;  0, 84, 327, 467;  0, 137, 172, 459;  0, 69, 88, 283;  0, 116, 251, 463;  0, 331, 421, 449;  0, 199, 207, 263;  0, 78, 301, 375;  0, 143, 221, 289;  0, 225, 295, 469;  0, 67, 111, 423;  0, 86, 177, 236;  0, 141, 158, 377;  0, 155, 345, 431;  0, 68, 257, 266;  0, 119, 217, 439}

{\boldmath $\adfPENT(4,156,13)$}, $d = 2$:
{115, 155, 170, 337;  115, 193, 309, 351;  115, 219, 319, 411;  115, 275, 312, 459;  155, 193, 319, 459;  155, 219, 312, 351;  155, 275, 309, 411;  170, 193, 312, 411;  170, 219, 309, 459;  170, 275, 319, 351;  193, 219, 275, 337;  309, 312, 319, 337;  337, 351, 411, 459;  24, 72, 132, 277;  24, 146, 207, 290;  24, 164, 264, 328;  24, 174, 208, 368;  72, 146, 264, 368;  72, 164, 208, 290;  72, 174, 207, 328;  132, 146, 208, 328;  132, 164, 207, 368;  132, 174, 264, 290;  146, 164, 174, 277;  207, 208, 264, 277;  277, 290, 328, 368;  192, 120, 431, 229;  351, 228, 174, 380;  158, 481, 268, 44;  369, 135, 442, 28;  21, 293, 160, 405;  52, 326, 406, 136;  470, 459, 296, 247;  37, 287, 76, 13;  61, 213, 378, 94;  95, 206, 429, 6;  467, 242, 305, 479;  105, 476, 140, 460;  473, 119, 380, 432;  415, 328, 230, 409;  46, 237, 278, 34;  2, 159, 301, 454;  146, 439, 126, 343;  133, 80, 205, 371;  382, 287, 477, 468;  410, 117, 34, 81;  173, 298, 118, 320;  327, 226, 176, 97;  65, 168, 292, 433;  254, 55, 464, 3;  241, 68, 302, 329;  386, 401, 132, 405;  256, 335, 333, 422;  0, 1, 6, 17;  0, 5, 35, 220;  0, 24, 172, 257;  0, 4, 134, 199;  0, 29, 51, 230;  0, 3, 201, 388;  0, 36, 271, 401;  0, 31, 251, 345;  0, 2, 279, 333;  0, 88, 209, 463;  0, 96, 383, 451;  0, 59, 143, 168;  0, 161, 227, 307;  0, 71, 247, 381;  0, 243, 263, 349;  0, 67, 369, 412;  0, 57, 66, 315;  0, 159, 176, 347;  0, 46, 405, 413;  0, 58, 385, 435;  0, 45, 188, 417;  0, 8, 215, 415;  0, 129, 397, 455;  0, 91, 112, 268;  0, 117, 429, 475}

{\boldmath $\adfPENT(4,160,13)$}, $d = 2$:
{39, 87, 137, 309;  39, 139, 247, 315;  39, 221, 293, 323;  39, 238, 256, 489;  87, 139, 293, 489;  87, 221, 256, 315;  87, 238, 247, 323;  137, 139, 256, 323;  137, 221, 247, 489;  137, 238, 293, 315;  139, 221, 238, 309;  247, 256, 293, 309;  309, 315, 323, 489;  6, 153, 172, 356;  6, 180, 248, 358;  6, 186, 343, 408;  6, 202, 274, 456;  153, 180, 343, 456;  153, 186, 274, 358;  153, 202, 248, 408;  172, 180, 274, 408;  172, 186, 248, 456;  172, 202, 343, 358;  180, 186, 202, 356;  248, 274, 343, 356;  356, 358, 408, 456;  260, 163, 223, 384;  151, 466, 55, 346;  337, 172, 149, 420;  347, 334, 283, 492;  391, 311, 44, 187;  264, 1, 441, 452;  161, 329, 438, 58;  415, 439, 50, 264;  251, 307, 5, 218;  139, 34, 35, 124;  209, 79, 4, 205;  466, 269, 424, 179;  385, 443, 432, 463;  28, 51, 310, 183;  240, 402, 442, 483;  319, 107, 16, 210;  365, 393, 136, 300;  206, 322, 263, 43;  386, 248, 295, 260;  468, 327, 489, 35;  228, 113, 103, 491;  206, 436, 225, 286;  26, 439, 492, 327;  346, 269, 154, 404;  4, 413, 100, 204;  185, 122, 52, 16;  298, 362, 493, 308;  38, 157, 355, 373;  0, 3, 24, 123;  0, 4, 220, 431;  0, 5, 364, 371;  0, 27, 56, 101;  0, 20, 187, 327;  0, 29, 140, 262;  0, 60, 241, 279;  0, 34, 345, 416;  0, 49, 74, 425;  0, 17, 148, 361;  0, 32, 159, 408;  0, 97, 304, 453;  0, 73, 305, 419;  0, 125, 161, 204;  0, 117, 129, 358;  0, 51, 173, 348;  0, 227, 421, 491;  0, 183, 225, 463;  0, 199, 357, 493;  0, 107, 342, 455;  0, 79, 111, 145;  0, 153, 198, 387;  0, 66, 321, 407;  0, 38, 275, 439;  0, 121, 128, 260;  0, 25, 269, 326}

{\boldmath $\adfPENT(4,164,13)$}, $d = 2$:
{67, 142, 213, 364;  67, 221, 271, 433;  67, 227, 309, 437;  67, 269, 281, 459;  142, 221, 309, 459;  142, 227, 281, 433;  142, 269, 271, 437;  213, 221, 281, 437;  213, 227, 271, 459;  213, 269, 309, 433;  221, 227, 269, 364;  271, 281, 309, 364;  364, 433, 437, 459;  48, 70, 74, 280;  48, 198, 238, 286;  48, 226, 263, 294;  48, 236, 245, 440;  70, 198, 263, 440;  70, 226, 245, 286;  70, 236, 238, 294;  74, 198, 245, 294;  74, 226, 238, 440;  74, 236, 263, 286;  198, 226, 236, 280;  238, 245, 263, 280;  280, 286, 294, 440;  407, 94, 371, 118;  449, 448, 345, 123;  95, 322, 424, 440;  374, 165, 285, 262;  493, 210, 317, 282;  35, 466, 101, 200;  499, 244, 381, 50;  350, 154, 241, 411;  114, 350, 298, 218;  231, 73, 354, 384;  227, 270, 236, 125;  241, 414, 335, 94;  218, 37, 167, 308;  1, 416, 415, 149;  169, 233, 405, 68;  481, 326, 500, 113;  486, 425, 242, 328;  145, 397, 83, 40;  289, 352, 444, 497;  471, 241, 268, 248;  94, 400, 279, 418;  457, 354, 452, 133;  401, 80, 278, 311;  131, 128, 386, 357;  81, 387, 338, 471;  109, 446, 93, 194;  137, 16, 215, 126;  0, 13, 78, 93;  0, 17, 130, 391;  0, 21, 41, 412;  0, 29, 32, 501;  0, 31, 358, 415;  0, 36, 81, 120;  0, 39, 113, 287;  0, 51, 245, 432;  0, 59, 83, 250;  0, 138, 401, 473;  0, 117, 276, 439;  0, 66, 143, 371;  0, 151, 349, 461;  0, 327, 379, 425;  0, 149, 208, 399;  0, 63, 207, 336;  0, 70, 176, 275;  0, 144, 301, 387;  0, 122, 323, 429;  0, 131, 265, 278;  0, 100, 259, 280;  0, 76, 192, 495;  0, 109, 134, 381;  0, 273, 343, 373;  0, 46, 218, 435;  0, 55, 108, 375;  1, 33, 109, 397;  0, 111, 272, 491;  0, 62, 126, 299}
\end{proof}
For the next three theorems, we provide a lemma concerning the existence of
group divisible designs with block size 4 and group type $g^u m^1$ when
$g$ and $m$ are even.
%
\begin{lemma}
\label{lem:4-GDD-existence}
Suppose $g$ and $m$ are even, $m \equiv g \adfmod{3}$ and $g \ge 8$.
Then there exists a {\rm 4-GDD} of type $g^{3t} m^1$ for
\begin{equation}
\begin{array}{lll}
t \ge 1  &\text{if} &m = g,\\
t \ge 4  &\text{if} &g < m \le g(3t-1)/2,\\
t \ge 10 &\text{if} &1 \le m < g.
\end{array}
\end{equation}
\end{lemma}
\begin{proof}
See \cite{BrouwerSchrijverHanani1977} or \cite[Theorem 4.6]{Ge2007} when $m = g$;
otherwise see \cite[Theorem 7.1]{WeiGe2015} and \cite[Theorem 1.2(i)]{ForbesForbes2018}.
\end{proof}
Although it is not needed for our purpose, if $g$ is divisible by 6,
we can improve Lemma~\ref{lem:4-GDD-existence} considerably; see \cite[Theorem 1.1]{ForbesForbes2018}.
%
\begin{theorem}
\label{thm:PENT-4-r-13-constructed-finite}
Let
$$R = \{112, 116, 120, 124, 128, 132, 136, 140, 144, 148, 152, 156, 160, 164\}.$$
Then there exists a $\adfPENT(4, (3r + 14) t + s, 13)$ with deficiency graph of girth $5$ for $r, s \in R$ and $t \ge t_0$,
where $t_0 = 1$ if $s = r$, $t_0 = 4$ if $s > r$, $t_0 = 10$ if $s < r$.
\end{theorem}
\begin{proof}
By Lemma~\ref{lem:PENT-4-r-13-direct}, there exists a $\adfPENT(4,r,13)$ with deficiency $(13,5)$-graph for each $r \in R$.

Suppose $r, s \in R$.
Take $3t$ copies of a $\adfPENT(4, r, 13)$, which has $3r + 14$ points, and a $\adfPENT(4,s,13)$, which has $3s + 14$ points.
By Lemma~\ref{lem:4-GDD-existence}, there exists a 4-GDD of type $(3r + 14)^{3t} (3s + 14)^1$ provided $t \ge t_0$,
where $t_0 = 1$ if $s = r$, $t_0 = 4$ if $s > r$, $t_0 = 10$ if $s < r$.
Therefore, by Theorem~\ref{thm:GDD-basic}, there exists a $\adfPENT(4, (3r + 14) t + s, 13)$ for $t \ge t_0$
with deficiency $(13,5)$-graph.
\end{proof}
%
\begin{theorem}
\label{thm:PENT-4-r-13-constructed}
There exists a $\adfPENT(4, r, 13)$ with deficiency graph of girth $5$ for all even $r \ge 253520$.
\end{theorem}
\begin{proof}
By Theorem~\ref{thm:PENT-4-r-13-constructed-finite} with $r = 116$ and $s = 120$,
there exists a $\adfPENT(4, 362t + 120, 13)$ with deficiency $(13,5)$-graph for $t \ge 4$.
Now take $3u$ copies of a $\adfPENT(4,112,13)$ from Lemma~\ref{lem:PENT-4-r-13-direct} and the $\adfPENT(4, 362t + 120, 13)$.
By Lemma~\ref{lem:4-GDD-existence}, there exists a 4-GDD of type $350^{3u} (1086t + 374)^1$ provided
$1086t + 374 \le 350(3u - 1)/2$, i.e.\ $u \ge (362t + 183)/175$.
Therefore, by Theorem~\ref{thm:GDD-basic},
there exists a $\adfPENT(4, 350u + 362t + 120, 13)$ with deficiency $(13,5)$-graph for $t \ge 4$ and $u \ge (362t + 183)/175$.
To show that all even $r \ge 253520$ are covered we argue as follows.

Dividing by two and ignoring the constant term, we consider
$$\{175u + 181t: t \ge 4, u \ge (362t + 183)/175\}.$$
Let $M = 175 \cdot 181$ and note that $\gcd(175, 181) = 1$.
By the Chinese Remainder Theorem, for any integer $m$, there exist $a_1$ and $a_2$ such that $m \equiv 175 a_1 + 181 a_2 \adfmod{M}$
with $0 \le a_1 < 181$ and $4 \le a_2 < 179$.
By adding a multiple of 181 to $a_1$ whenever necessary,
it follows that every integer $m$ in the range $M \le m < 2M$ can be expressed in the form
$m = 175 b_1 + 181 b_2$ with $b_1 \ge 0$ and $4 \le b_2 \le 178 = t_\mathrm{max}$.

Therefore every integer $n \ge 4M$ has a representation of the form $n = 175 (b_1 + 181h) + 181 b_2$, with
$b_1 \ge 0$, $4 \le b_2 \le t_\mathrm{max}$ and $h \ge 2$.
Put $t = b_2$ and $u = b_1 + 181h$.
Then $4 \le t \le t_\mathrm{max}$ and $u \ge 3 \cdot 181 > (362 t_\mathrm{max} + 183)/175$, as required.
Hence a $\adfPENT(4, r, 13)$ exists for all even $r \ge 8M + 120$.
\end{proof}
From the proof of Theorem~\ref{thm:PENT-4-r-13-constructed}
we see that the upper bound of admissible $r$ for which we fail to construct a $\adfPENT(4, r, 13)$ is rather large.
However, we can reduce the limit considerably if we drop the girth 5 condition and
permit the use of a $\adfPENT(4,4,13)$.
%
\begin{theorem}
\label{thm:PENT-4-r-13-constructed-degenerate}
There exists a $\adfPENT(4,r,13)$ if $r$ is even and $r \ge 474$.
\end{theorem}
\begin{proof}
Let $S = \{112, 116, 120, \dots, 160\}$ and note that $S$ covers the even residue classes modulo 26.
Let $t \ge 13$.

For each $s \in S$, take $3t$ copies of a $\adfPENT(4,4,13)$
and a $\adfPENT(4,s,13)$, which exists by Lemma~\ref{lem:PENT-4-r-13-direct}.
By Lemma~\ref{lem:4-GDD-existence}, there exists a 4-GDD of type $26^{3t} (3s + 14)^1$.
Hence, by Theorem~\ref{thm:GDD-basic} there exists a $\adfPENT(4, 26 t + s, 13)$.

The largest value of $26 t + s$ missed by this construction is $26 \cdot 12 + 160 = 472$.
The deficiency graphs of the constructed geometries consist of $3t$ components $K_{13,13}$ and
a 13-regular component with $3s + 14$ vertices and girth 5.
\end{proof}


\section{Concluding remarks}\label{sec:Concluding remarks}

Contrary to initial expectations, it has not been easy to find generalized pentagonal geometries
where the deficiency graph has girth at least 5 and the block size exceeds 3.
We have no examples of $\adfPENT(k,r,w)$ with $w > k \ge 5$, or with $w > 13$ and $k = 4$.
Even where we have been successful the current situation is not entirely satisfactory
in view of the considerable gap between
the smallest known $(w,5)$-graphs and the point set sizes of our geometries,
as is plainly illustrated by Table~\ref{tab:minimum-girth-5}.
\begin{table}[h]
\caption{Orders of $(w, 5^+)$-graphs}
\label{tab:minimum-girth-5}
\begin{tabular}{c|cccc}
     & Smallest & Moore & Smallest & Moore \\
     &  known   & bound &   known  & bound \\
 $w$ & $\adfPENT(3,r,w)$ & girth 6 & girth 5 & girth 5 \\
\hline
\end{document}